\documentclass[11pt,twoside]{article}
\usepackage{amsmath,amsfonts,amssymb,cases,mathdots,color,colordvi,url,tikz}
\usepackage{hyperref}
\usepackage{bm}
\usepackage{bbm}
\usepackage{pgfpages}
\usepackage{tabularx}
\usepackage{graphics}
\usepackage{subfigure}
\usepackage{epsfig}
\usepackage[figuresright]{rotating}
\usepackage{rotating}
\usepackage{algorithm}
\usepackage{algorithmic}
\usepackage{physics}
\usepackage{epstopdf}
\usepackage{parskip}
\usepackage{footmisc}
\usepackage{tablefootnote}
\usepackage[justification=centering]{caption}
\usepackage{booktabs}
\usepackage{amsmath}
\usepackage{array}
\usepackage{caption}
\usepackage{amsmath,amssymb}
\usepackage{booktabs}
\usepackage{graphicx}
\usepackage{tikz}

\usepackage{amsmath}
\usepackage{fix-cm}
\usepackage{xcolor}
\usepackage{hyperref}
\usepackage[capitalize,nameinlink,noabbrev]{cleveref}

\usetikzlibrary{arrows.meta}

\newcommand{\circled}[1]{%
	\tikz[baseline=(char.base)]{
		\node[shape=circle,draw,inner sep=1pt] (char) {#1};
	}%
}

\usepackage{amsthm}
\usepackage{aliascnt}
\usepackage{hyperref}
\usepackage[nameinlink,capitalize]{cleveref}

\newtheorem{proposition}{Proposition}

\newaliascnt{theorem}{proposition}
\newtheorem{theorem}[theorem]{Theorem}
\aliascntresetthe{theorem}

\newaliascnt{lemma}{proposition}
\newtheorem{lemma}[lemma]{Lemma}
\aliascntresetthe{lemma}

\newaliascnt{corollary}{proposition}

\aliascntresetthe{corollary}

\newaliascnt{definition}{proposition}
\newtheorem{definition}[definition]{Definition}
\aliascntresetthe{definition}

\newaliascnt{condition}{proposition}

\aliascntresetthe{condition}

\newaliascnt{conjecture}{proposition}

\aliascntresetthe{conjecture}

\newtheorem{remark}{Remark}
\newtheorem{assumption}{Assumption}

\crefname{proposition}{Proposition}{Propositions}
\Crefname{proposition}{Proposition}{Propositions}

\crefname{theorem}{Theorem}{Theorems}
\Crefname{theorem}{Theorem}{Theorems}

\crefname{lemma}{Lemma}{Lemmas}
\Crefname{lemma}{Lemma}{Lemmas}

\crefname{corollary}{Corollary}{Corollaries}
\Crefname{corollary}{Corollary}{Corollaries}

\crefname{definition}{Definition}{Definitions}
\Crefname{definition}{Definition}{Definitions}

\crefname{condition}{Condition}{Conditions}
\Crefname{condition}{Condition}{Conditions}

\crefname{conjecture}{Conjecture}{Conjectures}
\Crefname{conjecture}{Conjecture}{Conjectures}

\crefname{remark}{Remark}{Remarks}
\Crefname{remark}{Remark}{Remarks}

\crefname{assumption}{Assumption}{Assumptions}
\Crefname{assumption}{Assumption}{Assumptions}

\crefname{example}{Example}{Examples}
\Crefname{example}{Example}{Examples}

\newcommand{\be}{\begin{equation}}
	\newcommand{\ee}{\end{equation}}
\newcommand{\ba}{\begin{eqnarray}}
	\newcommand{\ea}{\end{eqnarray}}
\newcommand{\bas}{\begin{eqnarray*}}
	\newcommand{\eas}{\end{eqnarray*}}

\graphicspath{ {./Figures/} }

\def\T{{\mathbb T}}
\def\mbR{{\mathbb R}}
\def\R{{\mathcal R}}
\def\S{{\mathcal S}}

\def\D{{\mathcal D}}
\def\E{{\mathcal E}}

\def\I{{\mathcal I}}

\def\K{{\mathcal K}}
\def\cL{{\mathcal L}}

\def\N{{\mathcal N}}
\def\O{{\mathcal O}}

\def\T{{\mathcal T}}
\def\V{{\mathcal V}}

\def\bfp{{\bf p}}
\def\bfq{{\bf q}}
\def\bfr{{\bf r}}

\def\bfc{{\bf c}}
\def\bfd{{\bf d}}

\def\bft{{\bf t}}
\def\bfv{{\bf v}}
\def\bfx{{\bf x}}
\def\bfy{{\bf y}}
\def\bfz{{\bf z}}
\def\bfw{{\bf w}}
\def\bfu{{\bf u}}

\def\bfvsig{{\bf {\boldsymbol \varsigma}}}
\def\bfone{{\bf 1}}
\def\Prox{\mbox{Prox}}

\def\obz{\overline{\bfz}}
\def\bfdt{\boldsymbol{\delta}}

\def\dom{{\rm dom}\,}
\def\ri{{\rm ri}}

\def\wbz{\widetilde{\bfz}}
\def\wz{{\widetilde{z}}}
\def\whx{\widehat{\bfx}}
\def\why{\widehat{\bfy}}
\def\whz{\widehat{\bfz}}
\def\whw{\widehat{\bfw}}
\def\argmin{\mathop{{\rm argmin}}}

\def\whH{\widehat{H}}
\def\obx{\overline{\bfx}}
\def\oby{\overline{\bfy}}
\def\obz{\overline{\bfz}}
\def\obw{\overline{\bfw}}
\def\veps{\varepsilon}
\def\wtw{\widetilde{\bfw}}
\def\wtt{\widetilde{t}}
\def\wtx{\widetilde{\bfx}}
\def\dist{{\rm dist}}

\def\bp{ \textbf{Proof.} }
\def\ep{ \hfill $\Box$ }

\def\OT{ \overline{T} }

\def\mbN{ \mathbb{N} }

\def\obt{\overline{\beta}}

\def\bfla{\boldsymbol{\lambda}}
\def\bfxi{\boldsymbol{\xi}}
\def\bfom{\boldsymbol{\omega}}
\def\whp{\widehat{\partial}}
\def\bfvt{\boldsymbol{\vartheta}}
\def\bfmu{\boldsymbol{\mu}}

\def\pla{{\Phi_{\bfla}}}
\def\pmu{{\Phi_{\bfmu}}}
\def\prox{{\rm Prox}}

\def\tol{\texttt{tol}}
\def\mbone{\mathbbm{1}}

\title{A Proximal Point Semismooth Newton Algorithm for Structural Composite Cardinality Optimization}

\author{Penghe Zhang\thanks{Department of Data Science and Artificial Intelligence, The Hong Kong Polytechnic University, Hong Kong SAR, China, E-mail: {penghe.zhang@polyu.edu.hk} },
\ \ Naihua Xiu\thanks{School of Mathematics and Statistics, Beijing Jiaotong University, Beijing 100044, China, E-mail: {nhxiu@bjtu.edu.cn} } \ \ and \ \
Houduo Qi\thanks{Department of Data Science and Artificial Intelligence, and Department of Applied Mathematics, The Hong Kong Polytechnic University, Hong Kong SAR, China, E-mail: {houduo.qi@polyu.edu.hk} }
	}

\date{\today}

\begin{document}
	
	\maketitle	

	\begin{abstract}
Many image-processing problems can be formulated as composite cardinality optimization (CCOP), whose objective is the sum of two convex terms and a cardinality function composed with a linear operator. The composite cardinality term creates major algorithmic challenges: the separability of the cardinality function is lost and convergence analysis often requires surjectivity-type assumptions on the linear operators. To overcome these challenges, we consider the stationary dual formulation of CCOP, which has more favorable structure consisting of two convex terms and a separable cardinality function.
In this paper, we develop an efficient proximal point algorithm (PPA) to solve the stationary dual problem. The efficiency of our PPA stems from two aspects. Firstly, the key step of its subproblem solver minimizes a smooth convex function over a low-dimensional subspace by the classic semismooth Newton algorithm (SNA), which has global convergence and local superlinear rate under suitable conditions. Secondly, implementable inexact criteria are proposed for early termination of the SNA subroutine. These components form the basic framework of our inexact PPA. Under suitable conditions, it enjoys global convergence and local linear convergence rate.
In particular, we provide examples in which the convergence assumptions are automatically satisfied.
Finally, the SNA subroutine is incorporated into our inexact PPA to solve jump-sparse signal recovery and computed tomography (CT) image restoration. Numerical results demonstrate the time efficiency and solution accuracy of our proposed method.
	
	\vspace{3mm}
	
	\noindent{\bf \textbf{Keywords}:} nonconvex optimization, cardinaliry function, proximal point algorithm, semismooth Newton algorithm, inexact criteria, convergence analysis
	
\end{abstract}
{}

\section{Introduction} 
This paper focuses on efficient computation of the composite cardinality optimization (CCOP) as follows: 
\begin{align} \tag{P} \label{P}
	\min_{\bfx \in \mbR^n} f(\bfx) + g(A\bfx) + \pla(B\bfx), 
\end{align}
where $f:\mbR^n \to (-\infty, \infty]$ and $g:\mbR^m \to (-\infty, \infty]$ are proper, closed, and convex functions, $A: \mbR^n \to \mbR^m$ and $B : \mbR^n \to \mbR^p$ are linear operators, and $\pla:\mbR^p \to \mbR$ is the cardinality function defined by
\begin{equation}
	\pla(\bfu) := \sum_{i = 1}^p \lambda_i \mbone_{\{ u_i \neq 0 \}},
\end{equation}
where $\bfla : = (\lambda_1, \cdots, \lambda_p)^\top >0$ is the given regularization parameter vector,
and $\mbone_{\{ \O \}}$ is the characteristic function that takes 1 when the condition $\O$ holds and 0 otherwise.

Many practical applications arising from machine learning and image processing can be formulated as problem \eqref{P}. In particular, this model has flexible choices in imaging applications. The function $f$ serves as the convex regularizer or constraint, such as nonnegativity or box constraints on pixel intensities, to improve the performance of image processing. The function $g$ works as the loss function in the data-fidelity term. Its choice depends on the type of noise and the robust requirement of the model, including squared $\ell_2$, $\ell_1$, $\ell_\infty$, and Kullback-Leibler loss functions \cite{yuan2017ell,zhang2017nonconvex}. The linear operator $A$ determine the class of imaging tasks. For example, identity, convolution, and discrete Radon operators lead to denoising, deblurring, and CT reconstruction tasks respectively \cite{xu2011image,dong2013efficient,storath2015joint}. 
The composite term $\pla(B\bfx)$ promotes structural sparsity, where the choice of $B$ determines the relevant structure. For instance, finite-difference and wavelet-frame operators lead to sparsity of image gradients and frame coefficients respectively \cite{storath2014jump,zhang2013l0,dong2013efficient}.

Although problem \eqref{P} provides a unified framework for several imaging processing models, it is challenging to solve. The main difficulty comes from the composite cardinality term $\pla(B\bfx)$. The function $\pla$ is nonconvex and discontinuous, and its composition with a large-scale linear
operator $B$ brings substantial difficulties for both convergence analysis and numerical computation. In addition, the function $g$ is allowed to be
nonsmooth, and $f$ may represent additional convex constraints. This paper aims to overcome these difficulties and develop a globally convergent algorithm to find stationary solutions of problem \eqref{P}.



\subsection{Related works} \label{sec-related-works} In this section, we review the literature closely related to CCOP, including more general and special cases.



\textbf{Approximation approach.} In view of the difficulties arise in CCOP, many existing works adopt simplified models. For instance, the fidelity term $g$ is commonly specialized to the squared $\ell_2$ loss, while the term $f$, which may represent constraints such as box constraints, is sometimes omitted. Moreover, the cardinality term $\pla(B\bfx)$ is often replaced by convex or continuous approximation, such as $\ell_1$ total variation (TV) \cite{goldstein2009split}, $\ell_q$ TV with $0 < q < 1$ \cite{bian2015linearly,chen2012non}, truncated regularizer \cite{chouzenoux2013majorize,wu2018general}, and weighted difference regularizer \cite{lou2015weighted}. These approximate models have led to efficient algorithms and successful applications. Nevertheless, these models may not completely preserve the sparse structure of the original cardinality term. In particular, convex TV-type regularizer is known to exhibit contrast loss in some image restoration tasks \cite{wu2018general}.
This motivates the study of the original CCOP model.

\textbf{Decomposition approach.} A common approach for nonconvex composite optimization is to introduce auxiliary variables and enforce the resulting splitting constraints by Lagrangian or penalty functions. The Lagrangian-based methods includes the augmented Lagrangian method (ALM), the alternating direction method of multipliers (ADMM), and their variants. In recent years, this class of methods have been developed for problems involving the composition of a lower semicontinuous function with a linear or smooth mapping \cite{li2015global,bolte2018nonconvex,bot2019proximal,wang2019global,yashtini2021multi,song2020zero}. Their convergence analysis typically relies on the K$\L$ property, boundedness of primal and multiplier sequences, and regularity conditions on the data.
In particular, when applied to \eqref{P}, these methods usually require the linear operator $Q:\bfx\mapsto [A\bfx;B\bfx]$ to be surjective in order to establish whole-sequence convergence. 
A recent work \cite{hallak2023adaptive} bypass this surjectivity-type assumption by using an adaptive-feasibility mechanism to adjust the Lagrangian penalty parameter. However, their method only find an approximate critical point for each time.

Another related approach is the penalty-decomposition (PD) method. It is developed for cardinality-regularized nonlinear programming problems \cite{lu2013sparse} and has been successfully applied to image restoration models \cite{zhang2013l0}. PD generates a new iterate by applying block coordinate minimization to quadratic-penalty subproblems.
Similar to classic exterior penalty methods, PD enforces the splitting constraint by increasing the penalty parameter, which may be driven to infinity in the convergence analysis. 
Consequently, if high feasibility accuracy is required, subproblems become ill-conditioned. In addition, the existing convergence results for PD are mainly formulated in terms of accumulation points, rather than whole-sequence convergence. 

\subsection{Our approach and contribution}

\begin{figure}[h]
	\centering
	\includegraphics[width=1.0\textwidth]{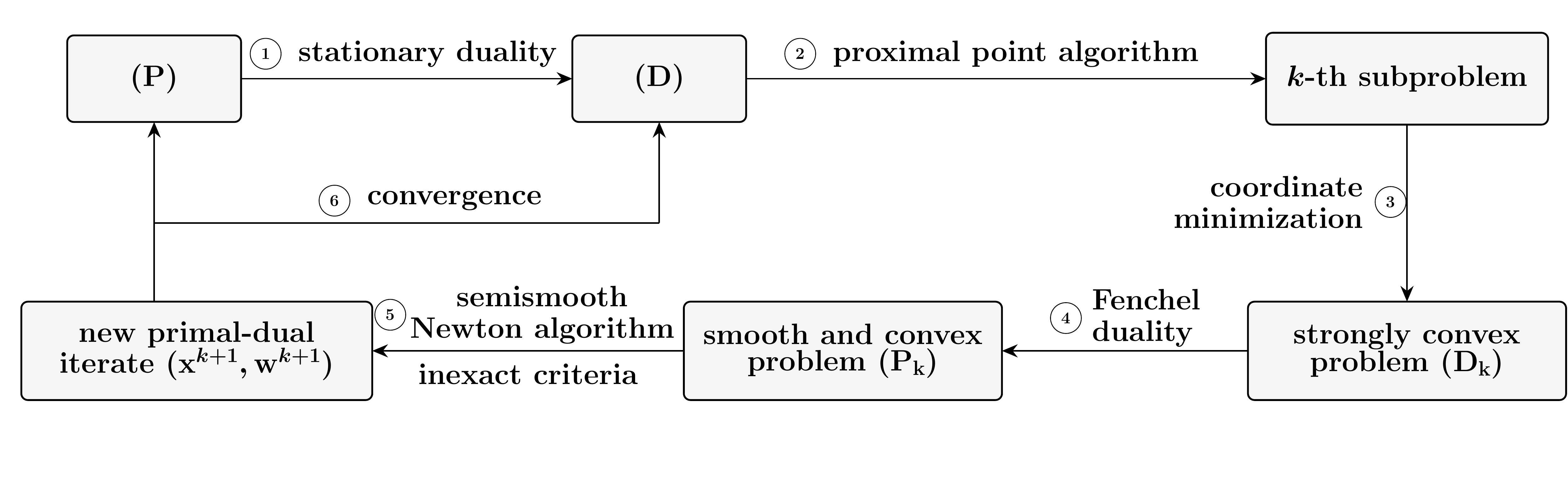}
	\caption{Flowchart of the proximal point algorithm in this work.}
	\label{flowchart}
\end{figure}

The literature review in the previous section indicates the methods for primal problem \eqref{P} typically rely on surjectivity-type regularity assumptions, boundedness of multipliers, or involve conditioning issues in the subproblems. These issues mainly arise from the composition of the cardinality function and a nontrivial linear operator. Alternatively, through the use of stationary duality theory in \cite{zhang2025composite, zhang2026local}, a recent work \cite{zhang2026stationaryduality} 
proposes a dual formulation of \eqref{P}:
\begin{align}  \tag{D} \label{D}
	\min_{\bfw = [\bfy; \bfz] \in \mbR^{m+p}} f^*(-Q^\top \bfw) + g^*(\bfy) + \pmu(\bfz),
\end{align}
where the linear operator $Q: \mbR^n \to \mbR^{m+p} $ is defined by $Q\bfx := [A\bfx; B\bfx]$, $f^*$ and $g^*$ are Fenchel conjugate functions of $f$ and $g$ respectively, and $\bfmu \in \mbR^p$ is a positive vector. It is noteworthy that the correspondence between local minimizers of \eqref{P} and \eqref{D} is established in \cite{zhang2026stationaryduality}. Moreover, the linear operator $B$ has been decoupled from the cardinality function $\pmu$, and thereby separable structure of $\pmu$ is available. 
These observations motivate us to solve the stationary dual problem \eqref{D}. 
However, the difficulty is that both $f^*$ and $g^*$ are unlikely differentiable. It is challenging to extract gradient or Hessian information from them and this prevents us from employing gradient-descent or Newton-type methods. 
To tackle those challenges,
we develop a general framework schemed in Fig.~\ref{flowchart}, which is briefly explained below.

{\bf (i) Efficient identification of the active subspace}.
The overall plan is to develop a proximal point algorithm (PPA) for \eqref{D}. 
Its subproblem at Step \circled{2} is
\begin{align} \tag{$k$-th subproblem} \label{kth-sub}
	\bfw^{k+1} \approx \argmin_{\bfw = [\bfy; \bfz]} \underbrace{f^*(-Q^\top \bfw) + g^*(\bfy) + \frac{1}{2\gamma_k} \|  \bfw - \bfw^k \|^2_M}_{\Xi_k(\bfw)} + \pmu(\bfz), 
\end{align}
where  $\bfw^k$ is the $k$th iterate, $\gamma_k > 0$, and  the linear operator $M: = QQ^\top + I$. Since $\Xi_k$ is nonsmooth and $\pmu$ is discontinuous, it is difficult to handle them simultaneously. In Step \circled{3}, by leveraging the separable and sparsity-induced structure of $\pmu$, we implement coordinate minimization on the subproblem and identify a low-dimensional subspace. 

{\bf (ii) Efficient solution of the subproblem over the subspace}.
We note that the function $\Xi_k$ restricted to the subspace is strongly convex and leads to a strongly convex problem \eqref{Dk}. The well-known Fenchel duality means that the dual problem of
\eqref{Dk} is convex and smooth \eqref{Pk} at Step \circled{4}.
For this class of convex programs, the classic semismooth Newton algorithm (SNA) enjoys global convergence with local superlinear rate. In particular, it has been successfully incorporated into PPA to solve large-scale sparse convex programming \cite{zhang2020proximal,yang2013proximal}.
Therefore, 
in Step \circled{5}, we use SNA to compute an approximate solution of \eqref{Pk} and take it as the new primal iterate $\bfx^{k+1}$. The associated dual iterate $\bfw^{k+1}$ can be recovered by applying the strong duality to convex problems \eqref{Pk} and \eqref{Dk}. 
The preceding steps describe how our PPA generate the primal-dual sequences $\{\bfx^k\}_{k \in \mbN}$ and $\{\bfw^k\}_{k \in \mbN}$. 

\textbf{(iii) Globally convergent PPA under inexact criteria}.
To further improve computational efficiency, we design verifiable inexact criteria for early termination of the SNA subroutine. 
This formulates our inexact PPA framework. Under a relative-interior-type regularity condition and other mild assumptions, we establish global convergence of our inexact PPA in Step \circled{6}: the dual sequence converges to a local minimizer of \eqref{D}, and every primal accumulation point is a local minimizer of \eqref{P}. If a local error bound is further satisfied, the dual sequence has local linear convergence rate.

\textbf{(iv) Weaker convergence assumptions}. We provide sufficient conditions and concrete examples under which the
relative-interior-type regularity condition and the error-bound condition
hold. In particular, these conditions hold automatically when \(f\) is the indicator
function of a box constraint containing the origin, and \(g\) is the
\(\ell_1\), squared \(\ell_2\), or \(\ell_\infty\) loss functions. In these cases, we do not
need the surjectivity-type regularity condition commonly imposed in
Lagrangian-based methods for nonconvex composite optimization.

Finally, the SNA subroutine is incorporated into the inexact PPA framework, resulting in a proximal point semismooth Newton algorithm (PPSNA) for CCOP. We further apply PPSNA to jump-sparse signal recovery and CT image restoration. Numerical comparisons with several leading solvers demonstrate the time efficiency and solution accuracy of PPSNA.

\subsection{Organization}
This paper is organized as follows. Some commonly used symbols, notations, and existing results in \cite{zhang2026stationaryduality} are briefly introduced in Section \ref{pre}. 
The implementation details of our PPA and convergence theory are given in Section \ref{sec-ppsna}. 
The numerical experiments on jump-sparse signal recovery and CT image restoration are conducted in Section \ref{sec-num-exp}. We conclude this paper in Section \ref{sec-conclu} with discussion on possible future research. For conciseness, all proofs in this paper are provided in the supplementary material.

\section{Preliminaries} \label{pre}
This section aims to describe the notation and definition frequently used in the paper. 

\subsection{Notation} \label{Subsection-Notation}
We denote $\mbR^n$ as $n$-dimensional Euclidean space endowed with an inner product $\langle \cdot, \cdot \rangle$ and the induced norm $\| \cdot \|$. The boldfaced lowercase letter $\bfx \in \mathbb{R}^n$ denotes a column vector of size $n$ and $\bfx^\top$ is its transpose. Let $x_i$ or $(\bfx)_i$ denote the $i$-th element of $\bfx$. Every positive semidefinite matrix $H$ induces a seminorm $\| \bfx \|_H:= \sqrt{\langle \bfx, H \bfx \rangle}$, where ``$:=$'' means ``define''. The $\delta$-neighborhood of $\bfx^* \in \mathbb{R}^n$ is denoted as $\mathcal{N}(\bfx^*, \delta) := \{ \bfx \in \mathbb{R}^n \ | \ \| \bfx - \bfx^* \| < \delta \}$. Given a set $\Omega$, we denote $\dist_H(\bfx,\Omega):= \inf_{\bfy \in \Omega}\| \bfx - \bfy \|_H $ and $\dist(\bfx,\Omega):= \inf_{\bfy \in \Omega}\| \bfx - \bfy \| $. The indicator function $\bfdt_\Omega(\bfx)$ takes 0 when $\bfx \in \Omega$ and $\infty$ otherwise. We use $\ri(\Omega)$ to represent the relative interior of $\Omega$ and its convex hull is $\mbox{co} (\Omega)$. We denote $\bfone$ as the vector with all elements being 1. The absolute value function is denoted as ${\rm abs}(\cdot)$.

Let $[p]$ be the set of indices $\{1, \ldots, p\}$.
For a subset $T \subset [p]$, $|T|$ denotes the number of elements of $T$, and $\OT$ consists of those indices of $[p]$ not in $T$. For vector $\bfz \in \mathbb{R}^p$ (resp. matrix $B \in \mbR^{p \times n}$), $\bfz_{T}$ (resp. $B_{T:}$) denotes the subvector of $\bfz$ indexed by $T$ (resp. the submatrix of $B$ with rows indexed by $T$). The largest and smallest eigenvalue function of a symmetric matrix are denoted by $\bfla_{\max}(\cdot)$ and $\bfla_{\min}(\cdot)$ respectively. 
The support set of $\bfz \in \mathbb{R}^p$ is denoted by $\mathcal{I}(\bfz):= \{ i\in [p]: z_i \neq 0  \}$. We denote $I$ as the identity matrix of appropriate dimension.

To simplify the notation, we denote $\bfw := [\bfy; \bfz] \in \mbR^{m+p}$ and define linear operator $Q: \mbR^n \to \mbR^{m+p} $ by $Q\bfx := [A\bfx; B\bfx]$. Moreover, we define 
\begin{align*}
	\Psi(\bfx) := f(\bfx) + g(A\bfx) ~~&\mbox{and}~~ \Xi(\bfw) := f^*(-Q^\top \bfw) + g^*(\bfy).
\end{align*} 
Then \eqref{P} and \eqref{D} can be represented as:
\begin{align*}
	\min_{\bfx } \Psi(\bfx) + \pla(B\bfx) ~~\mbox{and}~~ \min_{\bfw} \Xi(\bfw) + \pmu(\bfz). 
\end{align*}
%

\subsection{Variational tools} To proceed, we will introduce some definitions used throughout this paper.

\textbf{(A) Subdifferentials.} Let $\psi:\mbR^n\to(-\infty,\infty]$ be a proper lower semicontinuous (lsc) function, and define $\dom\psi:=\{\bfx\mid \psi(\bfx)<\infty\}$. For $\obx\in\dom\psi$, the limiting (or Mordukhovich) subdifferential \cite{mordukhovich2018variational,RockWets98} of $\psi$ at $\obx$ is denoted by $\partial\psi(\obx)$.
The limiting subdifferential of $\pmu$ has been calculated in \cite{le2013generalized}. Given $\bfz \in \mbR^p$ and $\bfmu > 0$, the representation is as follows:
\begin{align} \label{subdiff-l0}
	\partial \pmu(\bfz) = \left\{ \bfv ~\left|~ v_i \in \left\{ \begin{aligned}
		&  \mbR, && \mbox{if}~ z_i = 0, \\
		& 0,     && \mbox{otherwise},
	\end{aligned}  ~ i \in [p]  \right.  \right. \right\}. 
\end{align}
We can observe that $\partial \pmu$ is unchanged for any $\bfmu > 0$. 

\textbf{(B) Conjugacy.} Given $\psi:\mbR^n \to (-\infty,\infty]$, the conjugate function \cite{rockafellar1970convex} of $\psi$ is defined by
\begin{align} \label{conjugate-fun}
	\psi^*( \bfq ) := \sup_{\bfx \in \mbR^n} \langle \bfq, \bfx \rangle - \psi(\bfx).
\end{align}
One can verify the following Fenchel's inequality \cite[Theorem 4.6]{beck2017first} holds for any $\bfx, \bfq \in \mbR^n$:
\begin{align} \label{fenchel-ineq}
	\psi(\bfx) + \psi^*(\bfq) \geq \langle \bfx, \bfq \rangle.
\end{align}
If $\psi$ is proper, lsc, and convex, then $(\psi^*)^* = \psi$ and the following claims \cite[Theorem 23.5]{rockafellar1970convex} are equivalent for $\bfx,\bfq \in \mbR^n$:
\begin{align} \label{conjugate-subgradient}
	\psi(\bfx) + \psi^*(\bfq) = \langle \bfx, \bfq \rangle ~\iff~ \bfq \in \partial \psi(\bfx) ~\iff~ \bfx \in \partial \psi^*(\bfq).
\end{align}
Furthermore, given $\sigma > 0$, the conjugacy for $\sigma$-strong convexity and $(1/\sigma)$-smoothness (\cite[Definitions 5.1 and 5.16]{beck2017first}) has the following correspondence \cite[Theorem 5.26]{beck2017first}:
\begin{align*}
	\sigma\mbox{-strong convexity of}~\psi \iff (1/\sigma)\mbox{-smoothness of}~\psi^*.
\end{align*}

\textbf{(C) Moreau envelope.} For a proper and lsc function $\psi:\mbR^n \to  (-\infty,\infty]$, its Moreau envelope and proximal operator \cite{moreau1965proximite} at $\bfp \in \mbR^n$ with parameter $\sigma > 0$ are defined as follows:
\begin{align*}
	\E_{\sigma\psi}( \bfp ):= \min_{ \bfx } \frac{1}{2\sigma} \| \bfx - \bfp \|^2 + \psi( \bfx )~~\mbox{and}~~ \Prox_{\sigma\psi}(\bfp) := \argmin_{ \bfx } \frac{1}{2\sigma} \| \bfx - \bfp \|^2 + \psi( \bfx ).
\end{align*}
If $\psi$ is convex, then $\prox_{\sigma\psi}$ is single-valued and nonexpansive. Moreover, $\E_{\sigma\psi}$ is continuously differentiable \cite[Theorem 2.26]{RockWets98}
\begin{align} \label{moreau-gradient}
	\nabla \E_{\sigma\psi} (\bfp) = (1/\sigma) ( \bfp - \prox_{\sigma\psi} (\bfp) ).
\end{align}
The Moreau envelopes and proximal operators of $\psi$ and $\psi^*$ satisfies the following relationship \cite[Theorem 31.5]{rockafellar1970convex} in the convex setting:
\begin{align} \label{moreau-proximal-decomposition}
	\E_{\sigma\psi}(\bfp) + \E_{\frac{1}{\sigma} \psi^*} ( \bfp /\sigma ) =  \| \bfp \|^2/(2\sigma)~~\mbox{and}~~\prox_{\sigma\psi}(\bfp) + \sigma\prox_{\frac{1}{\sigma}\psi^*} (\bfp/\sigma) = \bfp
\end{align}

\textbf{(D) Semismoothness.} Let $\Omega$ be an open set and $F: \Omega \subseteq \mbR^n \to \mbR^n$ be locally Lipschitz continuous at $\bfx \in \Omega$. We denote $D_F$ as the set of differentiable points of $F$. Then the Clarke Jacobian \cite{clarke1990optimization} of $F$ at $\bfx$ is defined as $\partial_C F(\bfx) := {\rm co} \big\{ \lim_{\nu \to \infty} \nabla F( \bfx^\nu ) ~\big|~ \bfx^\nu \in D_F~\mbox{and}~ \bfx^\nu \to \bfx \big. \big\}.$
$F$ is said to be semismooth \cite{mifflin1977semismooth,qi1993nonsmooth} at $\bfx \in \Omega$ if $F$ is directionally differentiable at $\bfx$ and it holds that $F(\bfx + \bfd) - F(\bfx) - H\bfd = o(\| \bfd \|),~ \forall~ H \in \partial_C F(\bfx + \bfd).$

\subsection{Stationary dual theory of CCOP}
A recent work \cite{zhang2026stationaryduality} defines the stationary points of \eqref{P} and \eqref{D}, and further establish their correspondence. We summarize these useful results in this section. Firstly, let us present the definition of \eqref{P} and \eqref{D} respectively.
\begin{definition} \label{def-stationary-P}
	We say $\bfx^*$ is a stationary point of \eqref{P} if it satisfies 
	\begin{equation} \label{stationary-P}
		0 \in \partial f(\bfx^*) + A^\top \partial g(A\bfx^*) + B^\top \partial \pla(B \bfx^*).
	\end{equation}
\end{definition}
%
%
\begin{definition}
	We say $\bfw^*:= [\bfy^*;\bfz^*]$ is a stationary point of \eqref{D} if it satisfies
	\begin{align} \label{stationary-D}
		0 \in - \left( \begin{array}{c}
			A \\
			B
		\end{array} \right) \partial f^*( -Q^\top \bfw^* ) + \left( \begin{array}{c}
			\partial g^*(\bfy^*) \\
			\partial \pmu(\bfz^*)
		\end{array} \right). 
	\end{align}
\end{definition}

%
%
The relationship between stationary points and local minimizers of \eqref{P} and \eqref{D} are given in the following lemma.
\begin{lemma} \label{lem-correspondence}
	About problems \eqref{P} and \eqref{D}, we have
	
	(i) For each of the two problems, a stationary point implies a local minimizer.
	
	(ii) If $\bfx^*$ is a stationary point of \eqref{P}, then there exist $\bfy^* \in \partial g(A\bfx^*)$ and $\bfz^* \in \partial  \pla(B\bfx^*)$ such that $\bfw^* = [\bfy^*; \bfz^*]$ be a stationary point of \eqref{D}.
	
	(iii) If $\bfw^*$ is a stationary point of \eqref{D}, then there exists $\bfx^* \in \partial f^*(-Q^\top \bfw^*)$ such that $\bfx^*$ be a stationary point of \eqref{P}. 
\end{lemma}
Lem.~\ref{lem-correspondence} provides a theoretical guarantee for the validity of \eqref{D}. Particularly, a solution of \eqref{P} can be exactly recovered by solving \eqref{D} in certain cases. For example, when $f$ is strongly convex and a stationary point $\bfw^*$ of \eqref{D} is given, we can calculate a stationary point of \eqref{P} by $\bfx^* = \nabla f^*(-Q^\top \bfw^*)$. However, this becomes challenging when $f$ is merely convex and lsc. In the next section, we will design a PPA for \eqref{D}. Meanwhile, it generates a primal sequence whose accumulation points are stationary points of \eqref{P} under mild assumptions. 

\section{Proximal Point Semismooth Newton Algorithm} \label{sec-ppsna}
For the stationary dual problem \eqref{D}, all the linear terms are composited with convex function $f^*$, and the separable structure of $\pmu$ is available. 
%
%
%
%
%
%
%
%
%
%
%
Given the $k$-th dual iterate $\bfw^k:=[\bfy^k;\bfz^k]$, we consider developing a preconditioned proximal point algorithm (PPA) for \eqref{D}, which leads to the following problem
\begin{align} \label{Pk+card}
	\min_{\bfw = [\bfy;\bfz]}  \underbrace{\Xi(\bfw) + \frac{1}{2\gamma_k} \| \bfw - \bfw^k\|^2_M}_{:=\Xi_k(\bfw)} + \pmu(\bfz),
\end{align}
where $M:= Q Q^\top + I$ and therefore $\Xi_k$ is $(1/\gamma_k)$-strongly convex. The new iterates are generated by approximately solving \eqref{Pk+card}. However, this problem seems intractable due to the nonsmoothness of $\Xi_k$ and the discontinuity of $\pmu$. 

Fortunately, the cardinality term in \eqref{Pk+card} has separable structure and we can perform one cycle of coordinate minimization step on $\bfz$-variable for $i$ from 1 to $p$:
\begin{align} \label{coordinate-min}
	\wz^k_i \in \argmin\limits_{z_i \in \mbR}~ \underbrace{\Xi_k(\bfy^k, \wz^k_1, \cdots,\wz^k_{i-1}, z_i, z^k_{i+1},\cdots,z^k_p)}_{:=\bfxi_{ki}(z_i)} + \mu_i \mbone_{\{ z_i \neq 0 \}}.
\end{align}
The one-dimensional function $\bfxi_{ki}$ is strongly convex, and in many cases, we can find its global minimizer $\wtt^k_{i} = \argmin_{t} \bfxi_{ki}(t)$. Then one can verify that $\wz^k_i$ for $i \in [p]$ can be represented as
\begin{align} \label{sol-coordinate-min}
	\wz^k_i \in \left\{ \begin{aligned}
		& \wtt^k_{i}, && \mbox{if}~ \bfxi_{ki}(0) > \min \bfxi_{ki} + \mu_i, \\
		& \{ 0, \wtt^k_{i} \}, && \mbox{if}~ \bfxi_{ki}(0) = \min \bfxi_{ki} + \mu_i, \\
		& 0, && \mbox{otherwise}.
	\end{aligned} \right.
\end{align}
The derivation of \eqref{sol-coordinate-min} is given in the supplement materials. The iterate $\wbz^k = ( \wz^k_1, \cdots, \wz^k_p )^\top$ is used to help us identify an active set $T_k := \I(\wbz^k)$. After that, we renew the dual iterates by solving \eqref{Pk+card} restricted in the subspace $\bfz_{\OT_k} = 0$:
\begin{equation}
	\tag{D\textsubscript{k}}\label{Dk}
	\whw^{k+1}
	=
	\argmin_{\bfw = [\bfy;\bfz]}
	\Xi_k(\bfw)
	\quad \mbox{s.t.} \quad
	\bfz_{\OT_k} = 0 .
\end{equation}
Notably, \eqref{Dk} is strongly convex, and we only need to find an inexact solution $\bfw^{k+1} \approx \whw^{k+1}$ as a new dual iterate in practice. However, it is not an easy task to directly solve \eqref{Dk} because this problem involves composition of nonsmooth function and linear terms. Let us view \eqref{Dk} in an equivalent reformulation:
\begin{align*}
	\min_{\bfw = [\bfy;\bfz]}  f^*_k(-Q^\top \bfw) + g^*_k(\bfy) + h^*_k(\bfz), 
\end{align*} 
where
\begin{align*}
	& f^*_k(\bfp):= f^*(\bfp) + \frac{1}{2\gamma_k} \| \bfp + Q^\top \bfw^k \|^2,\\ 
	& g^*_k(\bfy) := g^*(\bfy) + \frac{1}{2\gamma_k} \|  \bfy - \bfy^k\|^2, \\
	& h^*_k(\bfz) := \frac{1}{2\gamma_k} \| \bfz - \bfz^k \|^2 + \bfdt_{\{0\}}(\bfz_{\OT_k}).
\end{align*}
As the three functions are all $(1/\gamma_k)$-strongly convex, 
the Fenchel dual problem of   \eqref{Dk} is smooth. We call it the  the primal problem corresponding to  \eqref{Dk}: 
\begin{align} \tag{P\textsubscript{k}} \label{Pk}
	\whx^{k+1} \in \argmin_{\bfx} \Psi_k(\bfx) := f_k(\bfx) + g_k(A\bfx) + h_k(B\bfx),
\end{align}
where by using the definition of conjugate function and property \eqref{moreau-proximal-decomposition}, we can compute
\begin{align*}
	& f_k(\bfx) := \E_{\frac{1}{\gamma_k} f} (\bfx - Q^\top \bfw^k/ \gamma_k) - \frac{1}{2\gamma_k} \| Q^\top \bfw^k \|^2, \\
	& g_k(\bfu):= \E_{\frac{1}{\gamma_k}g} (\bfu + \bfy^k/\gamma_k) - \frac{1}{2\gamma_k} \| \bfy^k \|^2, \\
	& h_k (\bfv):= \frac{\gamma_k}{2} \| (\bfv + \bfz^k/\gamma_k)_{T_k} \|^2 - \frac{1}{2\gamma_k}\| \bfz^k \|^2.
\end{align*}
All the three functions are $\gamma_k$-smooth, and particularly, $h_k$ is defined on a subspace $\mbR^{|T_k|}$. These structure allows us to apply a semismooth Newton algorithm (SNA) to solving \eqref{Pk} and \eqref{Dk} with fast convergence rate and low computational cost. Specifically, the subroutine SNA first generates a primal iterate by minimizing \eqref{Pk}, and then recovers an associated dual iterate by primal-dual solution correspondence between \eqref{Pk} and \eqref{Dk}.  

Theoretically, the sequence generated by SNA eventually converges to $(\whx^{k+1}, \whw^{k+1})$, but for numerical implementation, it returns an approximate iterates $(\bfx^{k+1}, \bfw^{k+1}) \approx (\whx^{k+1}, \whw^{k+1})$ as soon as certain stopping criteria are satisfied. It is worth mentioning that our PPA is implemented on a nonconvex and discontinuous problem \eqref{D}. Therefore, the stopping criteria for SNA are different from the convex setting \cite{rockafellar1976augmented,rockafellar1976monotone}. They are required to be more carefully designed so that our PPA has global convergence and local Q-linear convergence rate. We will give more detailed explanation on these criteria in Subsection \ref{sec-stopping-criteria}. Overall, the framework of PPA for \eqref{D} is summarized as follows.  
\begin{algorithm}[H] 
	\caption{PPSNA: a proximal point semismooth Newton algorithm for \eqref{D}.} \label{PPSNA}
	\begin{algorithmic}
		
		\STATE{Initialization: Input $(\bfx^0, \bfw^0)$ and $\gamma_0 > 0$.}
		\FOR{$k=0,1,\cdots$}
		
		\STATE{\textbf{1. Coordinate minimization: } } Compute $\wbz^k$ by \eqref{sol-coordinate-min} for $i \in [p]$ and select $T_k = \I(\wbz^k)$.
		

		\STATE{\textbf{2. Semismooth Newton:} } Denote $\wtw^k = [\bfy^k;\wbz^k]$ and apply the semismooth Newton algorithm (SNA, see Alg. \ref{SNA}) to \eqref{Pk} and obtain
		\begin{align*}
			(\bfx^{k+1}, \bfw^{k+1}) = \mbox{SNA}( \bfx^k, \wtw^k, \gamma_k, T_k )
		\end{align*}
		when certain stopping criteria are satisfied by $(\bfx^{k+1}, \bfw^{k+1})$ (see Subsection \ref{sec-stopping-criteria}).
		
		\STATE{\textbf{3. Parameter update :} }	Compute $\gamma_{k+1} \uparrow \gamma_\infty < \infty$.
		
		\ENDFOR
	\end{algorithmic}
\end{algorithm}

\subsection{Semismooth Newton algorithm for subproblems}
In this section, we make the following assumption.
\begin{assumption} \label{asm-coercive-f}
	Function $f$ is coercive (\cite[Definition 3.25]{RockWets98}), i.e. $\liminf_{\| \bfx \| \to \infty} f(\bfx)/\| \bfx \| = \infty$.
\end{assumption}
It implies coercivity of $\Psi_k$ by \cite[Exercise 3.29]{RockWets98}, and hence there must exist a global minimizer of \eqref{Pk}. We aim to design a semismooth Newton algorithm (SNA) which generates a sequence $\{ (\obx^j, \obw^j) \}_{j \in \mbN}$ converging to global minimizers $(\whx^{k+1}, \whw^{k+1})$ of \eqref{Pk} and \eqref{Dk}. Noticing that $f_k$, $g_k$, and $h_k$ are convex and $\gamma_k$-smooth, the optimality condition of \eqref{Pk} is given by
\begin{align} \label{equation-k}
	\nabla \Psi_k(\bfx) = \nabla f_k(\bfx) + A^\top \nabla g_k(A\bfx) + B^\top \nabla h_k(B\bfx) = 0, 
\end{align}
where by using properties \eqref{moreau-gradient} and \eqref{moreau-proximal-decomposition}, the gradient of $f_k$, $g_k$, and $h_k$ can be computed as follows: 
\begin{equation}\label{gradient-fghk}
	\begin{aligned} 
		& \nabla f_k (\bfx) = \gamma_k \Big( \bfx - Q^\top \bfw^k/\gamma_k - \prox_{\frac{1}{\gamma_k} f} (\bfx - Q^\top \bfw^k/\gamma_k) \Big) = \prox_{\gamma_k f^*} (\gamma_k \bfx - Q^\top \bfw^k), \\
		& \nabla g_k (\bfu ) = \gamma_k\Big(\bfu + \bfy^k/\gamma_k - \prox_{\frac{1}{\gamma_k}g} (\bfu + \bfy^k/\gamma_k)\Big)	= \prox_{\gamma_k g^*} (\gamma_k \bfu + \bfy^k), \\
		& \big(\nabla h_k (\bfv)\big)_{T_k} = (\gamma_k\bfv + \bfz^k)_{T_k} ~ \mbox{and}~ \big(\nabla h_k (\bfv)\big)_{\OT_k} = 0.
	\end{aligned}
\end{equation}
As $\nabla \Psi_k$ admits composite structure, it is difficult to give the exact form of Clarke Jacobian. We consider adopting the following generalized Jacobian:
\begin{align} \label{generalized-Jacob}
	\whp^2\Psi_k(\bfx) := \partial^2 f_k(\bfx) + A^\top \partial^2 g_k(A\bfx) A + \gamma_k B_{T_k:}^\top B_{T_k:}, 
\end{align}
where $\partial^2 f_k:= \partial_C(\nabla f_k)$ and $\partial^2 g_k:= \partial_C(\nabla g_k)$. Generally, $\partial^2\Psi_k$ and $\whp^2\Psi_k$ have the following relation (\cite[Page 75]{clarke1990optimization}):
\begin{align} \label{partial-whp}
	\partial^2\Psi_k(\bfx) \bfd \subseteq \whp^2\Psi_k(\bfx) \bfd,~~\forall \bfd \in \mbR^n.
\end{align}
From \eqref{equation-k}, \eqref{gradient-fghk} and \eqref{generalized-Jacob}, we have two important observations. On one hand, $\nabla \Psi_k$ is Lipschitz continuous, and also semismooth if both $\prox_{\gamma_k f^*}(\cdot)$ and $\prox_{\gamma_k g^*}(\cdot)$ are semismooth. On the other hand, since $\nabla_{\OT_k} h_k (\bfv) = 0$, the matrix $B$ can be replaced by the submatrix $B_{T_k:}$ when solving \eqref{equation-k}. These observations enable SNA to solve the subproblems with local superlinear convergence rate and lower computational cost. 

Let $\whx^{k+1}$ and $\whw^{k+1} = [\why^{k+1};\whz^{k+1}]$ be global minimizers of \eqref{Pk} and \eqref{Dk} respectively. Then according to \cite[Theorem 31.3]{rockafellar1970convex}, the following relationships hold: 
\begin{equation} \label{optimality-k}
	\begin{aligned} 
		&\why^{k+1} = \nabla g_k (A\whx^{k+1}),~ \whz^{k+1} = \nabla h_k(B\whx^{k+1}),~ \nabla f_k(\whx^{k+1}) + A^\top \why^{k+1} + B^\top\whz^{k+1} = 0, \\
		&G_k(\whx^{k+1}, \whw^{k+1}) := \Psi_k(\whx^{k+1}) + \Xi_k(\whw^{k+1}) = 0,
	\end{aligned}
\end{equation}
where the first line above reflects the stationarity of \eqref{Pk} and \eqref{Dk}, and $G_k:\mbR^{n}\times\mbR^{m+p} \to (-\infty,\infty]$ is the duality gap. As SNA generates a primal sequence $\{ \obx^j \}_{j \in \mbN}$ when solving \eqref{Pk}, we consider recover the corresponding dual iterates by \eqref{optimality-k} and \eqref{gradient-fghk}:	  
\begin{align} \label{dual-recovery}
	\oby^j  = \prox_{\gamma_k g^*} (\gamma_k A\obx^j + \bfy^k),~~ \obz^j_{T_k}  = \bfz^k_{T_k} + \gamma_k B_{T_k:} \obx^j ~~ \mbox{and}~~ \obz^j_{\OT_k} = 0. 
\end{align}
Now let us present the framework of SNA, which is adapted from \cite[Algorithm 2]{zhao2010newton}. 

		%
		%
		%
		%
		%
	
	\begin{algorithm}[H] 
		\caption{SNA: a semismooth Newton algorithm for subproblems (SNA($\bfx^k, \wtw^k, \gamma_k, T_k$)).} \label{SNA}
		\begin{algorithmic}
			
			\STATE{Initialization: Set starting point $(\obx^{0}, \obw^{0}) = (\bfx^k, \wtw^k)$, $\tau \in (0,1]$, $\tau_1, \tau_2, \beta \in (0,1)$, and $\zeta \in (0,1/2)$. }
			\FOR{$j=0,1,\cdots$}
			
			\STATE{\textbf{1. Newton direction: } } Compute $H_j \in \whp^2 \Psi_k (\obx^j)$ and $\theta_j =  \tau_1 \| \nabla \Psi_k(\obx^j) \|$. Apply the conjugate gradient method to find an approximate solution $\bfd^j$ to 
			\begin{align} \label{newton-equation}
				(H_j + \theta_j I) \bfd = - \nabla \Psi_k (\obx^j)
			\end{align}
			such that $\| (H_j + \theta_j I) \bfd^j + \nabla \Psi_k (\obx^j) \| \leq \tau_2 \| \nabla \Psi_k (\obx^j) \|^{1+\tau}$. 
			
			\STATE{\textbf{2. Line search:} } Set $\alpha_j = \beta^{q_j}$, where $q_j$ is the smallest nonnegative integer $q$ for which
			\begin{align*}
				\Psi_k ( \obx^j + \beta^q \bfd^j ) \leq \Psi_k(\obx^j) + \zeta \beta^q \langle \nabla \Psi_k(\obx^j), \bfd^j \rangle
			\end{align*}
			and set $\obx^{j+1} = \obx^j + \alpha_j \bfd^j$.
			
			\STATE{\textbf{3. Dual recovery:} }	Update dual iterate $\obw^j = [\oby^j; \obz^j]$ by \eqref{dual-recovery}. 
			
			\textbf{If} $(\obx^j,\obw^j)$ satisfies certain stopping criteria, \textbf{then} return $(\bfx^{k+1}, \bfw^{k+1}) = (\obx^j,\obw^j)$.
			
			
			
			\ENDFOR
		\end{algorithmic}
	\end{algorithm}
	
	\begin{remark} Here we give some comments about Alg. \ref{PPSNA}.
		
		(i) For step 1, We can notice that $H_j$ is positive semidefinite due to convexity of $f_k$ and $g_k$. Moreover, $\theta_j > 0$ if $\obx_j$ is not a solution to \eqref{Pk}. Taking these facts into account, Newton equation \eqref{newton-equation} is always solvable. For step 2, since $\Psi_k$ has Lipschitz continuous gradient, then according to \cite[Proposition 1.2.3]{bertsekas1997nonlinear}, we can find a uniform lower bound for $\{ \alpha_j \}_{j \in \mbN}$. Overall, SNA is well-defined.	
		
		(ii) Compared with \cite[Algorithm 2]{zhao2010newton}, our SNA adopts the generalized Jacobian \eqref{generalized-Jacob} tailored for \eqref{Pk}. By leveraging the sparsity of $\pmu$, one can observe that only subrows of $B$ are needed when computing \eqref{generalized-Jacob}. Moreover, SNA is equipped with a dual recovery step \eqref{dual-recovery}, and we will show that the dual sequence $\{ \obz^j \}_{j \in \mbN}$ converges to the unique global minimizer of \eqref{Dk}. 
	\end{remark}
	
	Next we show the global convergence of SNA.
	
	\begin{theorem} \label{thm-SNA-global-convergence}
		Suppose that Assumption \ref{asm-coercive-f} holds and let $\{(\obx^j, \obw^j)\}_{j \in \mbN}$ be the sequence generated by SNA, then the sequence is bounded and the following assertions hold:
		
		(i) Each accumulation point of $\{\obx^j\}_{j \in \mbN}$ is a global minimizer of \eqref{Pk}. We further have
		\begin{align} \label{Pk-convergence}
			\lim_{j \to \infty} \| \nabla \Psi_k (\obx^j) \| = 0 ~~ \mbox{and}~~ \lim_{j \to \infty} \Psi_k(\obx^j) = \min \Psi_k.
		\end{align}
		(ii) The sequence of duality gap $\{ G_k(\obx^j, \obw^j) \}_{j \in \mbN}$ converges to 0.
		
		(iii) $\{\obw^j\}_{j \in \mbN}$ converges to the unique global minimizer of \eqref{Dk} and we also have 
		\begin{align} \label{lim-Phik}
			\lim_{j \to \infty} \Xi_k(\obw^j) = \min\Xi_k.
		\end{align}
	\end{theorem}
	
	To show the local superlinear convergence of SNA, we need some assumptions of semismoothness and regularity.
	
	\begin{assumption} \label{asm-ssm-reg} The following conditions holds:
		
		(i) For any $\gamma >0$, both $\prox_{\gamma f^*}(\cdot)$ and $\prox_{\gamma g^*}(\cdot)$ are semismooth.
		
		(ii) Let $\{\obx^j\}_{j \in \mbN}$ be the primal sequence generated by SNA. There exists an accumulation point $\whx^{k+1}$ of the sequence such that every elements of $\whp^2 \Psi_k(\whx^{k+1})$ is nonsingular.
	\end{assumption}
	
	
	\begin{theorem} \label{thm-SNA-loc-superlinear}
		Suppose that Assumptions \ref{asm-coercive-f} and \ref{asm-ssm-reg} hold. Let $\{ (\obx^j, \obw^j) \}_{j \in \mbN}$ be the sequence generated by SNA. Then $\{ \obx^j \}_{j \in \mbN}$ converges to the global minimizer $\whx^{k+1}$ of \eqref{Pk} with Q-superlinear rate, i.e. 
		\begin{align} \label{SNA-superlinear}
			\| \obx^{j+1} - \whx^{k+1} \| = o(\| \obx^j - \whx^{k+1} \|).
		\end{align}
		Moreover, $\{ \obw^j \}_{j \in \mbN}$ converges to the global minimizer $\whw^{k+1}$ of \eqref{Dk} with R-superlinear rate.
	\end{theorem}
	
	\subsection{Implementable stopping criteria for the subroutine} \label{sec-stopping-criteria}
	In the last section, we show that SNA has global convergence and local superlinear convergence rate when solving subproblems. However, it eventually returns approximate solutions to \eqref{Pk} and \eqref{Dk} in practice. Therefore, the stopping criteria for the subroutines should be carefully designed to control the inexactness so that the global convergence of PPA can be achieved. Rockafellar proposed stopping criteria of subroutines for the classic PPA \cite{rockafellar1976monotone} and the analogy form for preconditioned PPA has been given in \cite{li2020asymptotically}. As for finding an approximate solution to \eqref{Dk}, the stopping criteria can be represented as:
	\begin{align}
		& \| \bfw^{k+1} - \whw^{k+1} \|_M \leq \veps_k,~0 \leq \veps_k,~\sum_{k=0}^{\infty} \veps_k < \infty, \label{C1-li}  \\
		& \| \bfw^{k+1} - \whw^{k+1} \|_M \leq \eta_k \| \bfw^{k+1} - \bfw^k \|_M,~ 0 \leq \eta_k < 1,~ \sum_{k=0}^{\infty} \eta_k < \infty, \label{C2-li}
	\end{align}
	where $\whw^{k+1}$ is the unique solution of \eqref{Dk}. These criteria are useful to analyze the convergence of PPA for convex programming, but they are difficult to be numerically implemented because $\whw^{k+1}$ is unknown in advance. Moreover, our problem \eqref{P} is nonconvex, the inexact criteria are required to be more carefully designed. Therefore, we consider adopting the following implementable criteria:
	\begin{align}
		& \Xi_k(\bfw^{k+1}) - \Xi_k(\wtw^k)  \leq \frac{1}{4\gamma_k} \| \bfw^{k+1} - \bfw^k \|^2_M, \tag{C1} \label{C1} \\
		&	G_k( \bfx^{k+1}, \bfw^{k+1} )  \leq \frac{\veps_k^2}{2\gamma_k}, ~0 \leq \veps_k,~\sum_{k=0}^{\infty} \veps_k < \infty, \tag{C2} \label{C2} \\
		&		G_k( \bfx^{k+1}, \bfw^{k+1} )  \leq \frac{\eta_k^2}{2\gamma_k} \| \bfw^{k+1} - \bfw^k \|^2_M,~ 0 \leq \eta_k < 1,~ \sum_{k=0}^{\infty} \eta_k < \infty, \tag{C3}  \label{C3}
	\end{align}
	We can observe that \eqref{C1} requires that the objective value of \eqref{Dk} at semismooth Newton iterate can be controlled. \eqref{C2} and \eqref{C3} guarantee the duality gap is small enough. Particularly, these two conditions imply \eqref{C1-li} and \eqref{C2-li} respectively. This can be derived by the strong convexity of $\Xi_k$ and weak duality between \eqref{Pk} and \eqref{Dk}:
	\begin{align*}
		G_k(\bfx^{k+1}, \bfw^{k+1}) = \Xi_k(\bfw^{k+1}) + \Psi_k(\bfx^{k+1}) \geq \Xi_k(\bfw^{k+1}) - \min\Xi_k  \geq \frac{1}{2\gamma_k} \| \bfw^{k+1} - \whw^{k+1} \|^2_M.
	\end{align*}
	Moreover, through the use of Thm. \ref{thm-SNA-global-convergence}, we can prove that the iterates generated by SNA satisfies the three inexact criterion after finite iterations.
	\begin{theorem} \label{thm-finite-termination}
		Suppose that Assumption \ref{asm-coercive-f} holds and let $\{(\obx^j, \obw^j)\}_{j \in \mbN}$ be the sequence generated by SNA. Then for sufficiently large $j$, $(\bfx^{k+1}, \bfw^{k+1}) = (\obx^j, \obw^j)$ satisfies \eqref{C1}, \eqref{C2}, and \eqref{C3}.
	\end{theorem}
	
	
	\subsection{Convergence analysis of PPSNA} 
	In this subsection, we aim to analyze the global convergence and local linear convergence rate of PPSNA when the subroutine SNA adopt stopping criteria \eqref{C1}, \eqref{C2} or/and \eqref{C3}. We need to stress that our PPSNA is for a nonconvex and discontinuous problem \eqref{D}. Even though the subproblems \eqref{Pk} and \eqref{Dk} are convex, the index set $T_k$ may keep changing, which brings significant challenges on the convergence analysis. Moreover, the inexact criteria \eqref{C1}--\eqref{C3} are designed specifically for our SNA subroutine, and hence differ from those in the existing literature \cite{attouch2009convergence,attouch2013convergence}. As a result, there is no off-the-shelf proof. For the convenience of readers, we list the main procedure of our proof.
	
	\textbf{Step 1.} The dual objective value sequence is decreasing and convergent (see Lem. \ref{lem-c1}).
	
	\textbf{Step 2.} For sufficiently large $k$, set $T_k$ remains unchanged  (see Lem. \ref{lem-c1-c2}). As a result, PPSNA becomes a preconditioned PPA for a convex program \eqref{D-T-infty}.
	
	\textbf{Step 3.} The convergence of preconditioned PPA \cite{li2020asymptotically} for convex programming is available to our PPSNA. Its global convergence with local linear rate is derived in Thms. \ref{thm-global-convergence} and \ref{thm-convergence-rate}.

	In our analysis, we require that Assumption \ref{asm-coercive-f} holds. On one hand, it ensures SNA has global convergence and satisfies the stopping criteria after finite iteration. On the other hand, this assumption implies $\dom f^* = \mbR^n$ from \cite[Theorem 11.8]{RockWets98}, which guarantees that each dual iterate $\bfw^k$ is feasible for \eqref{D}. For the conciseness of proof, we make the following assumption, and verifiable version will be provided in Remark \ref{rem-assumption}. 
	%
	%
	%
	\begin{assumption} \label{asm-bound}
		The following conditions hold.
		
		(i) Sequence $\{ \bfw^k \}_{k \in \mbN}$ is bounded.  
		
		(ii) Function $g^*$ is continuous on its domain, which means that for any $\oby$ and $\{ \bfy^\nu \}$ belonging to $\dom g^*$ with $\bfy^\nu \to \oby$, it holds that  $\lim_{\bfy^\nu \to \oby} g^*(\bfy^\nu) = g(\oby)$.
	\end{assumption}
	%
	
	\begin{lemma} \label{lem-c1}
		Suppose that Assumptions \ref{asm-coercive-f} and \ref{asm-bound} hold. If the stopping criteria \eqref{C1} is executed in Alg. \ref{SNA}, then the following assertions hold for Alg. \ref{PPSNA}: 
		
		(i) The dual objective value sequence is convergent and it holds that
		\begin{align*} 
			&\Xi(\bfw^k) + \pmu(\bfz^k) \geq  \Xi(\bfw^{k+1}) + \pmu(\bfz^{k+1}) + \frac{1}{4\gamma_k} \| \bfw^{k+1} - \bfw^k \|^2_M ~\mbox{and}~ \lim_{k \to \infty} \| \bfw^{k+1} - \bfw^k \| =0.
		\end{align*}	
		(ii) $\pmu(\bfz^{k})$ remains unchanged for all sufficiently large $k$. 
		%
	\end{lemma}
	
	When \eqref{C1} and \eqref{C2} are used as stopping criteria of SNA, we have the following result.
	
	\begin{lemma} \label{lem-c1-c2}
		Suppose that Assumptions \ref{asm-coercive-f} and \ref{asm-bound} hold. If the stopping criteria \eqref{C1} and \eqref{C2} are adopted by Alg. \ref{SNA}, then the following assertions hold for Alg. \ref{PPSNA}:
		
		(i) It holds that
		\begin{align} \label{wtw-wk}
			\lim_{k \to \infty} \| \wtw^{k} - \bfw^k \| =0.
		\end{align}		
		(ii) There exists a constant $\vartheta_1 >0$ such that $| \wz^k_i | \geq \vartheta_1$ for every $i \in T_k$ and $k \in \mbN$.
		
		(iii) When $k$ is sufficiently large, there exists a fixed index set $T_\infty \subseteq [p]$ such that
		\begin{align*}
			T_k = T_\infty.
		\end{align*}
	\end{lemma}
	
	Lems. \ref{lem-c1} and \ref{lem-c1-c2} indicate that Alg. \ref{PPSNA} always implies sufficiently descent on objective function value of \eqref{D}, then after finite iteration, $T_k$ is fixed and Alg. \ref{PPSNA} eventually becomes a preconditioned PPA for the following convex program:
	\begin{align} \label{D-T-infty}
		\min_{\bfw = [\bfy;\bfz]}  \Xi(\bfw)~ s.t.~ \bfz_{\OT_\infty} = 0.
	\end{align}
	If we regard it as a dual formulation, the corresponding primal program is as follows: 
	\begin{align} \label{P-T-infty}
		\min_{\bfx}  \Psi(\bfx)~ s.t.~ (B\bfx)_{T_\infty} = 0.
	\end{align}
	Now we give the following assumption.
	\begin{assumption} \label{asm-DT-solvable}
		The set $\{ \bfx \ | \ \bfx \in \ri(\dom f),~A\bfx \in \ri(\dom g),~ (B\bfx)_{T_\infty} = 0 \}$  is nonempty. Particularly, ``$\ri$" can be omitted if the corresponding $f$ or $g$ is polyhedral convex.
	\end{assumption}
	If Assumptions \ref{asm-coercive-f} and \ref{asm-DT-solvable} hold, then according to \cite[Corollary 31.2.1]{rockafellar1970convex}, the optimal value of \eqref{D-T-infty} is finite and attainable. Since $\dom f^* = \mbR^n$, the optimality condition of \eqref{D-T-infty} is as follows:
	\begin{align*}
		0 \in \S(\bfw):= - \left( \begin{array}{c}
			A \\
			B
		\end{array} \right) \partial f^*( -Q^\top \bfw ) + \left( \begin{array}{c}
			\partial g^*(\bfy) \\
			\V_\infty
		\end{array} \right), 
	\end{align*}
	where $\V_\infty = \{ \bfv \in \mbR^p \ | \ \bfv_{T_\infty} = 0 \}$. In other words, the solution set \eqref{D-T-infty} is $\S^{-1}(0)$. Next let us show the global convergence of PPSNA.
	\begin{theorem} \label{thm-global-convergence}
		Suppose that Assumptions \ref{asm-coercive-f}, \ref{asm-bound}, and \ref{asm-DT-solvable} hold. Moreover, criteria \eqref{C1} and \eqref{C2} are adopted by Alg. \ref{SNA}. Let $\{(\bfx^{k}, \bfw^{k})\}_{k \in \mbN}$ be the sequence generated by Alg. \ref{PPSNA}, then we have
		
		(i) $\{\bfw^{k}\}_{k \in \mbN}$ converges to a stationary point of \eqref{D}, which is also a local minimizer.
		
		(ii) Each accumulation point of $\{\bfx^{k}\}_{k \in \mbN}$ is a stationary point of \eqref{P}, which is also a local minimizer.
	\end{theorem}
	
	To establish the local Q-linear convergence rate of PPSNA, let us give the following assumption proposed in \cite{li2020asymptotically}.
	\begin{assumption} \label{asm-error-bound}
		For any $\nu_1 >0$, there exists $\zeta_1 > 0$ such that 
		\begin{align*}
			\dist( \bfw, \S^{-1} (0) ) \leq \zeta_1 \dist( 0, \S(\bfw) ), ~~ \mbox{whenever}~ \dist(\bfw, \S^{-1} (0)) \leq \nu_1.
		\end{align*}
	\end{assumption}
	Li et al. \cite{li2020asymptotically} showed that a sufficient condition to the above assumption is the local upper Lipschitz continuity \cite{robinson1976implicit} of $\S^{-1}$ at the origin. Robinson \cite{robinson2009some} established the celebrated result that every polyhedral multifunction is locally upper Lipschitz continuous. The graph of polyhedral multifunction are unions of finitely many polyhedral convex sets. This class of multifunctions has extensive applications and are closed under finite addition, scalar multiplication, and finite composition. Taking these facts into account, we give the following proposition.
	\begin{proposition} \label{prop-error-bound}
		If $\partial f^*$ and $\partial g^*$ are polyhedral multifunctions, then Assumption \ref{asm-error-bound} holds.
	\end{proposition}
	
	\begin{theorem} \label{thm-convergence-rate}
		Suppose that Assumptions \ref{asm-coercive-f}, \ref{asm-bound}, \ref{asm-DT-solvable}, and \ref{asm-error-bound} hold. Moreover, criteria \eqref{C1}, \eqref{C2}, and \eqref{C3} are adopted by Alg. \ref{SNA}. Let $\{\bfw^{k}\}_{k \in \mbN}$ be the dual sequence generated by Alg. \ref{PPSNA}, then for sufficiently large $k$, we have
		\begin{align*}
			\dist_M( \bfw^{k+1}, \S^{-1}(0) ) \leq \rho_k \dist_M( \bfw^k, \S^{-1}(0) ),
		\end{align*}
		where $\rho_k = (1 - \eta_k)^{-1} \big( \eta_k + (1+\eta_k) \zeta_1 \bfla_{\max}(M)/\sqrt{\gamma_k^2 + \zeta_1^2 \bfla^2_{\max}(M)} \big) $ and 
		\begin{align} \label{convergence-factor}
			\limsup_{k \to \infty} \rho_k = \frac{\zeta_1 \bfla_{\max}(M)}{\sqrt{\gamma_\infty^2 + \zeta_1^2 \bfla^2_{\max}(M)}} < 1.
		\end{align}
	\end{theorem}
	
	
	At the end of this section, let us provide easy-to-check assumptions and examples that guarantee the convergence of our PPSNA.
	
	\begin{remark} \label{rem-assumption}
		According to the proofs in our supplement materials, Assumption \ref{asm-bound} is mainly used to derive the following facts:
		\begin{itemize}
			\item[$(a)$] $\Xi$ is bounded below.
			\item[$(b)$] $f^*$ is Lipschitz continuous on the set containing all iterates generated by the coordinate minimization and semismooth Newton steps in PPSNA.
			\item[$(c)$] $\{ \bfy^k \}_{k \in \mbN}$ is bounded and $g^*$ is continuous on its domain. 
		\end{itemize}
		Here we provide the following sufficient conditions for (a), (b), and (c) respectively.	
		\begin{itemize}
			\item[$(a')$] $\{ \bfx \ | \ \bfx \in \ri(\dom f),~A\bfx \in \ri(\dom g),~ B\bfx = 0 \} \neq \emptyset,$ where ``$\ri$" can be omitted if the corresponding $f$ or $g$ is polyhedral convex.
			\item[$(b')$] $f^*$ is Lipschitz continuous.
			\item[$(c')$] $g^*$ is level-coercive and continuous on its domain. 
		\end{itemize}
		It follows from \cite[Assumption 2, Theorem 4]{zhang2026stationaryduality} that $(a')$ implies $(a)$. Moreover, it is also a sufficient condition of Assumption \ref{asm-DT-solvable}. As for $(c')$, the level coercivity of $g^*$ means that $g^*(\bfy) \to \infty$ if $\| \bfy \| \to \infty$. 
		
		Taking the above facts into consideration, we can use Assumption \ref{asm-coercive-f}, $(a')$, $(b')$, and $(c')$ to derive the global convergence of PPSNA in Thm. \ref{thm-global-convergence}. If we further assume $\partial f^*$ and $\partial g^*$ are polyhedral multifunctions, then the dual sequence generated by PPSNA has local linear convergence rate according to Thm. \ref{thm-convergence-rate}. In particular, these conditions automatically holds when $f$ is the indicator function of a box constraint containing the origin point, and $g^*$ is taken as $\ell_1$, $\ell_\infty$, or squared $\ell_2$ loss. Notably, these examples require neither surjectivity-type assumptions on the involved linear operators nor boundedness of primal and multiplier sequences, which are commonly used in Lagrangian-based methods for nonconvex optimization.

		
	\end{remark}

	\section{Numerical Experiments} \label{sec-num-exp}
	In this section, we give more details on numerical implementation of PPSNA. To evaluate the performance of our PPSNA, we compare it with other algorithms on solving jump-sparse signal recovery and CT image restoration problems. All the experiments are implemented on Matlab 2023b by a laptop with 32GB memory and Intel CORE i7 2.6 GHz CPU. 
	
	\subsection{Warm Start of PPSNA}
	In practice, warm-starting helps PPA to achieve fast local convergence earlier, which is crucial to its numerical efficiency. A common strategy is to first apply a simple first-order method, typically ADMM, to produce a moderate-accuracy initial point, and then invoke PPA to achieve faster local convergence and attain high accuracy efficiently. 
	Such a two-stage design is common in recent numerical implement for large-scale convex optimization problems \cite{li2020asymptotically,zhang2020proximal}.
	
	To implement a proximal ADMM (PADMM) for solving \eqref{P}, we consider equivalently reformulating it as the following constrained optimization
	\begin{align*}
		\min_{\bfx, \bfq, \bfu, \bfv} f(\bfq) + g(\bfu) + \pla(\bfv),~s.t. \ \bfq = \bfx, \ \bfu = A \bfx, \  \bfv = B\bfx.
	\end{align*}
	The corresponding augmented Lagrangian function can be represented as:
	\begin{align*}
		\cL_\varrho(\bfx, \bfq, \bfu, \bfv, \bfp, \bfy, \bfz): = & f(\bfq) + g(\bfu) + \pla(\bfv) + \langle \bfp, \bfx - \bfq \rangle + \langle \bfy, A\bfx - \bfu \rangle + \langle \bfz, B\bfx - \bfv \rangle \\
		&+ \frac{\varrho}{2} \| \bfx - \bfq \|^2 + \frac{\varrho}{2} \| A\bfx - \bfu \|^2+ \frac{\varrho}{2} \| B\bfx - \bfu \|^2 
	\end{align*}
	For convenience, we denote the $l$-th iterate as $\bfom^{l}: = (\bfx^l,\bfq^l,\bfu^l, \bfv^l, \bfp^l, \bfy^l, \bfz^l)$. Then the framework of our proximal ADMM is presented in Alg. \ref{PADMM}.
	\begin{algorithm}[htbp] 
		\caption{PADMM: a proximal ADMM for warm-start of PPSNA} \label{PADMM}
		\begin{algorithmic}	
			\STATE{Initialization: Input $\bfom^0$, $\obt$, and $\varrho_0 > 0$.}
			\FOR{$l=0,1,\cdots$}
			\STATE{\textbf{1. Primal step: } } Compute $\bfv^{l+1}$ and $(\bfx^{l+1}, \bfq^{l+1}, \bfu^{l+1})$.
			\begin{align}
				\bfv^{l+1} \in & \min_\bfv \cL_{\varrho_l}(\bfx^l,\bfq^l,\bfu^l, \bfv, \bfp^l, \bfy^l, \bfz^l) \label{v-step} \\
				(\bfx^{l+1}, \bfq^{l+1}, \bfu^{l+1}) \in & \min_{\bfx, \bfq, \bfu} \cL_{\varrho_l} ( \bfx, \bfq, \bfu, \bfv^{l+1}, \bfp^l, \bfy^l, \bfz^l) + \frac{\sigma}{2} \| \bfx - \bfx^l \|^2 \label{x-q-u-step}
			\end{align}	
			\STATE{\textbf{2. Multiplier step:} } Compute $\bfp^{l+1}$, $\bfy^{l+1}$, and $\bfz^{l+1}$.
			\begin{align*}
				\bfp^{l+1} = \bfp^l + \varrho_l(\bfx^{l+1} - \bfq^{l+1}),  \bfy^{l+1} = \bfy^l + \varrho_l (A\bfx^{l+1} - \bfu^{l+1}),  \bfz^{l+1} = \bfz^l + \varrho_l(B\bfx^{l+1} - \bfv^{l+1}).	 
			\end{align*}	
			\STATE{\textbf{3. Parameter update :}}	Compute $\Delta^{l+1}:= \| \bfom^{l+1} - \bfom^l \|/(1+ \|\bfom^l\|)$ and update $\varrho_{l+1}$ by
			\begin{align*}
				\varrho_{l+1} = \left\{\begin{aligned}
					& \min\{10^3, \obt\varrho_{l} \}, && \mbox{if } l \geq 1 \mbox{ and } \Delta^{l+1} > \Delta^{l} \\
					& \varrho_{l}, && \mbox{otherwise}.
				\end{aligned} \right.
			\end{align*}
			\ENDFOR
		\end{algorithmic}
	\end{algorithm}
	By using the proximal operator and Moreau envelop, \eqref{v-step} can be computed by 
	\begin{align*}
		\bfv^{l+1} \in & \Prox_{\frac{1}{\varrho_l} \pla} (B \bfx^{l+1} + \bfz^l/\varrho_l)  
	\end{align*}
	and \eqref{x-q-u-step} can be simplified as
	\begin{align*}
		&\bfx^{l+1} \in \argmin_\bfx \E_{\frac{1}{\varrho_l} f} ( \bfx + \bfp^l/ \varrho_l ) + \E_{\frac{1}{\varrho_l} g} ( A\bfx + \bfy^l/ \varrho_l ) + \frac{\varrho_l}{2} \| B\bfx + \bfz^l/\varrho_l - \bfv^{l+1} \|^2 + \frac{\sigma}{2} \| \bfx - \bfx^l \|^2,  \\
		& \bfq^{l+1} = \Prox_{\frac{1}{\varrho_l} f} (\bfx^{l+1} + \bfp^l/\varrho_l) \ \mbox{and} \ \bfu^{l+1} = \Prox_{\frac{1}{\varrho_l} g} (A\bfx^{l+1} + \bfy^l/\varrho_l).
	\end{align*}
	Particularly, $\bfx^{l+1}$ can be efficiently computed by SNA (see Alg. \ref{SNA}). As we introduced in Section \ref{sec-related-works}, the global convergence of ADMM in nonconvex setting may require surjective-type assumption on the linear operators or sequence boundedness \cite{bolte2018nonconvex,bot2019proximal,wang2019global}. These conditions may affect its practical performance. However, we emphasize that our PADMM is only used to generate an initial point to warm-start PPSNA. Therefore, our PADMM is not expected to achieve global convergence. It is also noteworthy that we increase the penalty parameter $\varrho_l$ dynamically, which is different from the classic ADMMs. This scheme helps to adjust the sparsity term $\pla(B\bfx^{l})$ according to our empirical tests. The detailed parameter setting of our PADMM will be given in the next subsections.  

	\subsection{Jump-sparse signal recovery} \label{exp-jump-sparse}
	Jump-sparse signals arise in many signal-processing tasks where the ground-truth signal is approximately piecewise constant and contains only a small number of abrupt transitions. Typical examples include stepwise biological and imaging data \cite{shepp1974fourier}. In practice, the ground-truth signal is often required to be recovered from blurred, incomplete, and noisy observation \cite{storath2014jump}. Specifically, given the ground-truth $\bfx^*$ in the range $[0,1]$, the observation $\bfc$ is obtained by
	\begin{align} \label{c-signal-recovery}
		\bfc = \T \bfx^* + \bfvt,
	\end{align} 
	where $\T \in \mbR^{m \times n}$ is a reduced blur matrix associated with the blur kernel $\bft:= (t_{-r}, \cdots,t_0, \cdots, t_r)$, defined by
	\begin{align} \label{A-signal-recovery}
		\T_{i,j} =  \left\{ \begin{aligned}
			& t_{i-j}, && \mbox{if}~ |i - j| \leq r, \\
			& 0, && \mbox{else},
		\end{aligned} \right.
	\end{align}  
	and $\bfvt$ is Laplacian noise with independent zero-mean components and variance $\upsilon^2$. To recover the underlying piecewise-constant signal, we consider solving the following problem
	\begin{align} \label{jump-sparse-opt}
		\min_{\bfx \in \mbR^n} \| \T\bfx - \bfc\|_1 + \pla(D_n\bfx), \quad s.t. \ 0 \leq \bfx \leq \bfone,
	\end{align}
	where $D_n \in \mbR^{(n-1) \times n}$ is the one-dimensional finite difference operator, defined by
	\begin{align} \label{Dn}
		D_n =
		\begin{bmatrix}
			-1 & 1 & 0 & \cdots & 0 \\
			0 & -1 & 1 & \ddots & \vdots \\
			\vdots & \ddots & \ddots & \ddots & 0 \\
			0 & \cdots & 0 & -1 & 1
		\end{bmatrix}
	\end{align}
	We can observe that \eqref{jump-sparse-opt} is a special case of \eqref{P} with the setting of $f(\bfx) = \bfdt_{[0,1]^n} (\bfx)$, $g(\bfu) = \| \bfu - \bfc \|_1$, $A = \T$, and $B = D_n$.
	
	Next, let us give the parameter setting of algorithms. For warm-start stage, we set $\obt = 1.05$, $\sigma = 0.1$, and $\varrho_0 = 1$. Denoting $\bfom^{l}: = (\bfx^l,\bfq^l,\bfu^l, \bfv^l, \bfp^l, \bfy^l, \bfz^l)$, the stopping criteria in this stage is as follows:
	\begin{align*}
		\frac{\| \bfom^{l+1}  - \bfom^l \| }{1 + \| \bfom^{l} \|} < \texttt{tol}_1
	\end{align*}
	We set a moderate tolerance $\tol_1 = 10^{-2}$ for PADMM. Now let $\bfom^l$ be a warm-start solution computed by PADMM. For the second stage to solve \eqref{D}, we take $\bfmu = \beta_1 \texttt{abs}(\bfv^l) + \beta_2$, where \texttt{abs} takes absolute value for each component of the given vector, $\beta_1 = 1$, and $\beta_2 =  \min \{|(B\bfx)_i| /(pn) \ | \ (B\bfx)_i \neq 0 \}$. About PPSNA, we set $\gamma_0 = 1$ and then $\gamma_k$ is updated by a dynamic strategy. At each iteration $k$, the update period is set to 2, 3, 5, and 10 for $k < 10$, $10 \leq k < 20$, $20 \leq k < 100$, and $k \geq 100$ respectively. Whenever $k$ is divisible by the given update period, we compute $\gamma_{k} = \min\{ 2 \gamma_{k-1}, 10^7 \}$, otherwise $\gamma_k$ is unchanged. For the inexact criteria \eqref{C2} and \eqref{C3}, we set $\varepsilon_k = \eta_k = 1/1.05^k$. Moreover, to implement SNA, we take $\tau_1 = 10^{-6}$, $\tau_2 = 0.1$, $\beta = 0.5$, and $\vartheta = 10^{-4}$. For problem \eqref{D}, we define the violation of dual optimality as follows
	\begin{align*}
		{\rm VDO}:=\max\left\{
		R_1^{k+1},\, R_2^{k+1},\, R_3^{k+1}
		\right\},
	\end{align*}
	where
	\begin{align*}
		R_1^{k+1}
		&:=
		\frac{
			\left\|
			Q^\top \bfw^{k+1}
			+
			\Prox_{f^*}\bigl(-Q^\top \bfw^{k+1}+\bfx^{k+1}\bigr)
			\right\|
		}{
			1+\left\|Q^\top \bfw^{k+1}\right\|
		},                                                    \\[0.5em]
		R_2^{k+1}
		&:=
		\frac{
			\left\|
			\bfy^{k+1}
			-
			\Prox_{g^*}\bigl(\bfy^{k+1}+A\bfx^{k+1}\bigr)
			\right\|
		}{
			1+\left\|\bfy^{k+1}\right\|
		},~ ~R_3^{k+1}
		:=
		\frac{
			\left\|(B\bfx^{k+1})_{T_k}\right\|
		}{
			1+\left\|\bfz^{k+1}\right\|
		}.                                                 	
	\end{align*}
	We stop PPSNA when ${\rm VDO} < \tol_2$ and our default tolerance is $\tol_2 = 10^{-4}$. As for the coordinate descent step in PPSNA, we give the closed-form formula and complexity analysis in the supplement materials.
	
	Now let us introduce the comparison methods. Pottslab \cite{storath2014jump} is a Matlab/Java toolbox for reconstructing jump-sparse signals and images based on the inverse Potts model. In particular, a problem similar to \eqref{jump-sparse-opt} is studied in \cite{storath2014jump}, but without the box constraint. CVX \cite{grant2014cvx} is a Matlab-based software package that can efficiently solve a broad class of convex optimization problems through interfaces to second-order interior-point based solvers such as SDPT3 \cite{toh2012implementation} and SeDuMi \cite{sturm1999using}. PD3O \cite{yan2018new} is a first-order primal dual algorithm for minimizing the sum of three functions with a linear operator. In particular, we apply CVX and PD3O to solve the convex optimization which replaces $\pla(D_n\bfx)$ in \eqref{jump-sparse-opt} by $\ell_1$ loss with weight vector $\bfla$.

	%
	
	As for the data in experiments of this section, the kernel $\bft$ is generated by first drawing a random vector $\bfvsig : = ( \varsigma_{-r}, \cdots, \varsigma_0, \cdots, \varsigma_r ) $ with entries uniformly sampled from $[0,1]$. The kernel entries are computed by $t_i = \exp(-\varsigma_i^2/2)$ for $i = -r, \cdots, 0, \cdots, r$, and then normalized so that the sum of all entries equals to one. Particularly, we set $r = 20$ in our experiments. The ground-truth signal $\bfx^*$ is constructed as a one-dimensional piecewise-constant signal with $n_J$ jumps. Specifically, $n_J + 1$ segment heights are randomly sampled from $[0,1]$. The jump locations are then uniformly selected from the points $1, \cdots, n$. Finally, we generate matrix $A$ and vector $\bfc$ by formulas \eqref{A-signal-recovery} and \eqref{c-signal-recovery} respectively. 
	\subsubsection{Illustration on warm start and convergence} \label{illustration-convergence}
\begin{figure}[htbp]
		\centering
		\begin{minipage}[t]{0.47\linewidth}
			\centering
			\includegraphics[width=\linewidth]{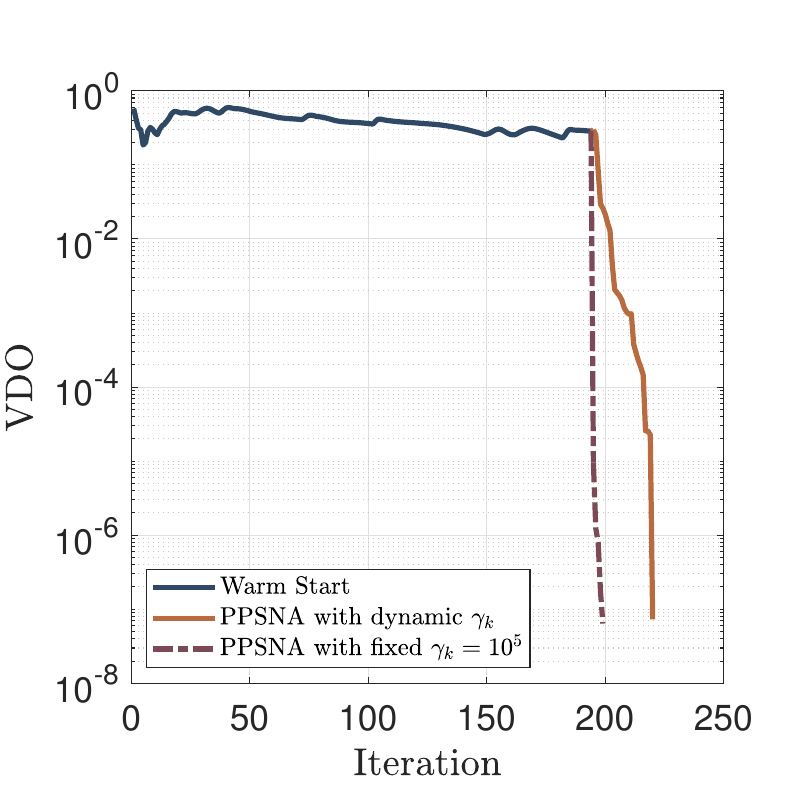}
		\end{minipage}
		\hfill
		\begin{minipage}[t]{0.47\linewidth}
			\centering
			\includegraphics[width=\linewidth]{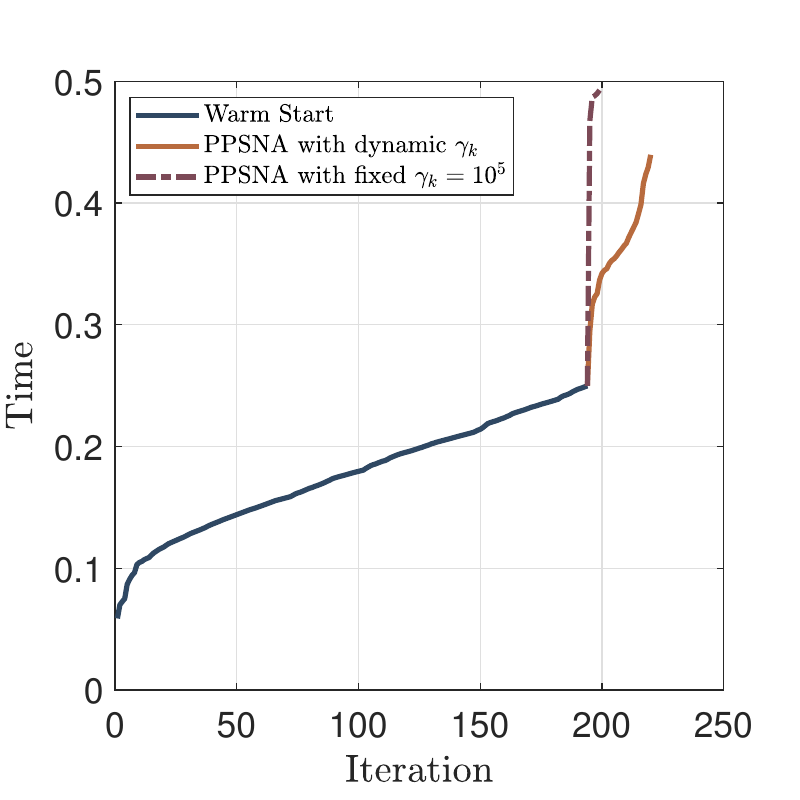}
		\end{minipage}
		\caption{Change of VDO and Time along with iteration of PPSNA on the data with
			$m=100$, $n=200$, $\upsilon=0.05$, and $n_J=5$.}
		\label{fig-convergence}
	\end{figure}

	There are two goals in this part. First, we will demonstrate the convergence rate of PPSNA at different stages. For this, we record the changes of VDO (resp. computational time) along with iteration of PPSNA. Second, we will explore how the dynamic setting of $\gamma_k$ influence the convergence rate. For this, we also fixed $\gamma_k = 10^5$ and then make comparison with the two strategies. To conduct these tests, we generate the data $A$ and $\bfc$ with $m = 100$, $n = 200$, $\upsilon = 0.05$, and $n_J = 5$. In this experiment, we take $\bfla = \lambda \bfone$ with $\lambda = 0.03$. The results are presented in Fig. \ref{fig-convergence}, and we denote the computational time as Time in the following experiments.
	
	As we can see from the first panel of Fig. \ref{fig-convergence}. The warm-start stage of PPSNA only has slight decrease on VDO, where as the second stage contribute to the significant decrease of VDO to achieve the satisfactory tolerance. Moreover, when fixed a large $\gamma_k = 10^5$, the convergence of PPSNA becomes much faster because according to Thm. \ref{thm-convergence-rate}, the Q-linear convergence factor in \eqref{convergence-factor} is smaller, leading to a faster convergence rate along with iteration. 
	
	However, this setting is more time-consuming than the dynamic strategy. The reason is that a large $\gamma_k$ often leads to a large condition number of \eqref{newton-equation}. As a result, PPSNA require more Time to solve the subproblem at each iteration. In the following experiments, we therefore adopt the dynamic strategy, which gradually increases $\gamma_k$ to reduce Time.
	
	
	%
	
	\subsubsection{Comparisons on data with various dimensions}
\begin{figure}[htbp]
		\centering
		\includegraphics[width=0.9\textwidth]{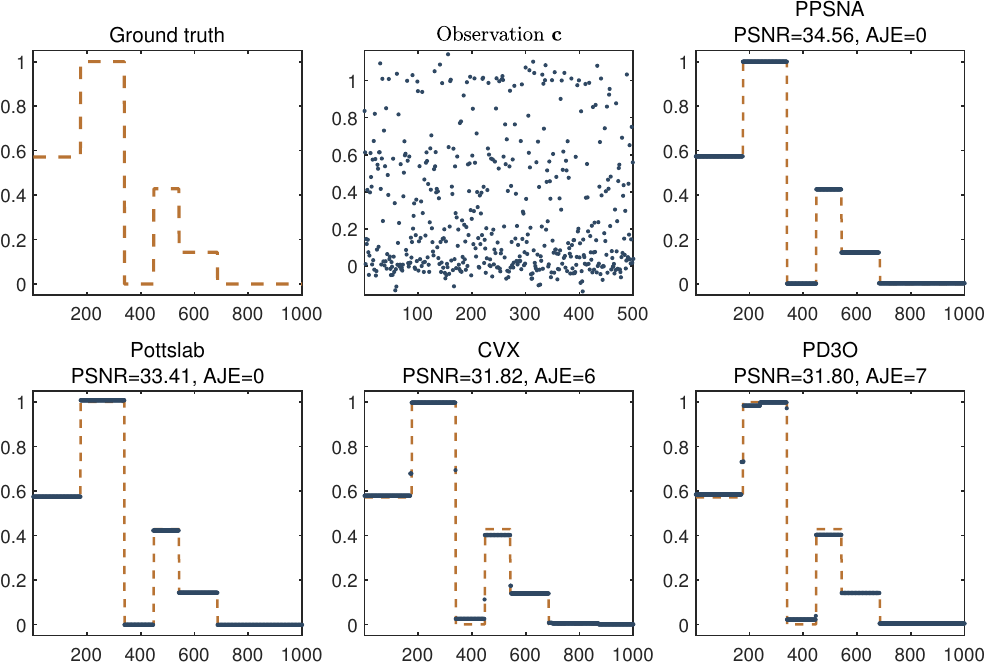}
		\caption{PSNR and AJE of all algorithms on data with $m = 500$, $n = 1000$, $\upsilon = 0.05$, and $n_J = 5$.}
		\label{fig-PC-signals}
	\end{figure}

	In this part, we will compare PPSNA with other algorithms on data with different $m$, $n$, and $\upsilon$. In this experiment, we take $\bfla = \lambda \bfone$ with $\lambda > 0$. To achieve better performance of all the algorithms, their regularization parameter $\lambda$ are all selected from $\{ 10^{-2}, 2\times 10^{-2}, \cdots, 10^{-1}, 2 \times 10^{-1}, \cdots, 1, 2, \cdots, 10 \}$.
	Besides computational time (Time), we record peak signal-to-noise-ratios (PSNR) at $\bfx$:
	\begin{align*}
		{\rm PSNR}(\bfx) := 10 \log_{10} \left( n \frac{\| \bfx^* \|^2_\infty}{\| \bfx - \bfx^* \|^2} \right),
	\end{align*} 
	where $\bfx^*$ is the ground-truth. The absolute jump-count error (AJE) is computed by
	\begin{align*}
		{\rm AJE} (\bfx) := {\rm abs} \left( \sum_{i = 1}^p \mbone_{\{ (B\bfx)_i \neq 0 \}} - \sum_{i = 1}^p \mbone_{\{ (B\bfx^*)_i \neq 0 \}} \right).
	\end{align*}
	We can notice that a larger PSNR indicates that the recovered signal $\bfx$ is closer to the ground truth. A lower AJE means that the recovered signal better captures the piecewise structure of the ground-truth signal. In the process of experiments, we notice the solutions computed by CVX and PD3O have significant incorrect jumps. Therefore, we refine them by the following procedure. First, given the solutions $\bfx^{sol}$, we identify the index set $T_{sol} := \{ i \in [n-1] ~|~ | (B \bfx^{sol})_i | < 10^{-3} \}$, which has sufficiently small finite differences. Then we compute the following projection to obtain the refined solution $\obx$:
	\begin{align*}
		\obx := \argmin_{\bfx} \| \bfx - \bfx^{sol} \| \quad s.t.~ (B\bfx)_{T_{sol}} = 0.
	\end{align*}
	After this process, the refined signals of CVX and PD3O also have piecewise-constant structure with much lower AJE.
	To illustrate PSNR and AJE, we first conduct experiments on a small-scale data with $m = 500$, $n = 1000$, $\upsilon = 0.05$, and $n_J = 5$. The results are demonstrated in Fig. \ref{fig-PC-signals}.  
	We can observe that ${\rm AJE} = 0$ for PPSNA and Pottalab. Therefore, compared with CVX and PD3O, their solutions achieve more accurate reconstructions of the ground truth, with larger PSNR.
	
	Next let us conduct the numerical comparison on data with larger $m$, $n$, and $\upsilon$. The results are summarized in Tab.~\ref{tab-dimension}. Now we give some comments about the results. For $m \in \{ 3000, 4000, 5000 \}$, $n = 10000$, and $\upsilon = 0.05$, PPSNA achieves the highest PSNR and zero AJE, while also exhibiting very competitive runtime, comparable to CVX and noticeably lower than PD3O and Pottslab. For $n \in \{2000,4000,6000\}$, $m = 2000$, and $\upsilon = 0.05$, PPSNA again delivers the best PSNR with zero AJE for all tested problem. Its runtime is almost 1/5 of Pottslab. For $\upsilon \in \{ 0.03,0.05,0.07 \}$, $m = 1500$, and $n = 3000$, PPSNA remains robust in both reconstruction quality and jump recovery when the noise level increase.
	
	\begin{sidewaystable}[htbp]
		\caption{Numerical results of all algorithms on data with different $m$, $n$, and $\upsilon$.}
		\label{tab-dimension}
		\centering
		\setlength{\tabcolsep}{2.8pt}
		\renewcommand{\arraystretch}{1.05}
		\begin{tabular*}{\textwidth}{@{\extracolsep{\fill}}ccccccccccccc@{}}
			\toprule
			& \multicolumn{4}{c}{PSNR $\uparrow$} & \multicolumn{4}{c}{AJE $\downarrow$} & \multicolumn{4}{c}{Time (sec) $\downarrow$} \\
			\cmidrule(lr){2-5} \cmidrule(lr){6-9} \cmidrule(lr){10-13}
			& PPSNA & Pottslab & CVX & PD3O & PPSNA & Pottslab & CVX & PD3O & PPSNA & Pottslab & CVX & PD3O \\
			\midrule
			\multicolumn{1}{c}{$m$} & \multicolumn{12}{c}{$n = 10000$, $\upsilon = 0.05$}  \\
			\midrule
			3000 & 34.83 & 33.74 & 34.51 & 33.83 & 0 & 0 & 25 & 8 & 1.732e+1 & 1.003e+3 & 1.697e+1 & 2.221e+1 \\
			4000 & 35.16 & 33.62 & 35.15 & 34.27 & 0 & 0 & 25 & 6 & 2.037e+1 & 2.516e+3 & 2.314e+1 & 3.216e+1 \\
			5000 & 36.40 & 34.40 & 35.50 & 34.05 & 0 & 0 & 25 & 5 & 4.301e+1 & 3.014e+3 & 3.051e+1 & 4.845e+1 \\
			\addlinespace[0.25em]
			\midrule
			\multicolumn{1}{c}{$n$} & \multicolumn{12}{c}{$m = 2000$, $\upsilon = 0.05$} \\
			\midrule
			2000 & 30.72 & 28.30 & 29.11 & 28.94 & 0 & 0 & 13 & 7 & 4.749e+0 & 7.875e+2 & 6.076e+0 & 1.998e+1 \\
			4000 & 31.71 & 31.64 & 31.52 & 29.93 & 0 & 0 & 18 & 5 & 9.287e+0 & 4.970e+2 & 9.178e+0 & 2.262e+1 \\
			6000 & 32.08 & 31.53 & 31.77 & 31.75 & 0 & 0 & 23 & 7 & 8.855e+0 & 4.718e+2 & 9.569e+0 & 1.710e+1 \\
			\addlinespace[0.25em]
			\midrule
			\multicolumn{1}{c}{$\upsilon$} & \multicolumn{12}{c}{$m = 1500$, $n = 3000$} \\
			\midrule
			0.03 & 32.15 & 32.96 & 32.45 & 32.46 & 0 & 0 & 18 & 6 & 6.078e+0 & 3.299e+2 & 7.202e+0 & 1.421e+1 \\
			0.05 & 30.78 & 30.57 & 30.33 & 29.97 & 0 & 0 & 17 & 6 & 5.441e+0 & 3.313e+2 & 6.414e+0 & 1.810e+1 \\
			0.07 & 28.68 & 28.91 & 28.59 & 26.65 & 0 & 0 & 15 & 3 & 5.824e+0 & 3.572e+2 & 6.449e+0 & 2.397e+1 \\
			\bottomrule
		\end{tabular*}
		{\\ \par\vspace{0.1cm} ``$\uparrow$" means that a larger metric value corresponds to better algorithm performance, while ``$\downarrow$" has the opposite meaning.}
	\end{sidewaystable}
	
	%
	%
	%
	%
	\subsection{CT image restoration with sparse angles}
	Computed tomography (CT) reconstruction from sparse-angle projections is a representative limited-data inverse problem. Practically, sparse-angle projections are often used to reduce radiation exposure. However, incompleteness of the projection data makes the reconstruction severely ill-posed. For this reason, it is natural to incorporate structural priors into the reconstruction model, especially when the underlying image is expected to be approximately piecewise constant. We consider the following sparse-angle CT restoration model \cite{storath2015joint}
	\begin{align} \label{ct-restoration}
		\min_{\bfx \in \mbR^n} \frac{1}{2} \| \R \bfx - \bfc \|^2 + \pla(\D \bfx), ~s.t.~ 0 \leq \bfx \leq 255.
	\end{align}
	Here, $\bfx \in \mbR^n$ denotes the vectorized CT image with $n = n_1n_2$ being the number of pixels in an $n_1 \times n_2$ image, $\R \in \mbR^{m \times n}$ is the sampled Radon transform with total number of measurement $m = Ls$, where $L$ is the number of detector bins and $s$ is the number of projection angles. Moreover, $\bfc \in \mbR^m$ is the observed projection data, and $\D$ is the two-dimensional finite-difference operator defined by
	\begin{align}
		\D := \left( \begin{array}{c}
			I_{n_2} \otimes D_{n_1} \\
			D_{n_2} \otimes I_{n_1},
		\end{array} \right)
	\end{align}
	where $I_{n_1}$ denotes the $n_1 \times n_1$ identity matrix and $D_{n_1}$ has been defined in \eqref{Dn}. The constraint $0 \leq \bfx \leq 255$ reflects the 8-bit grayscale range of the image data. We can observe that \eqref{ct-restoration} is a special case of \eqref{P} with the setting of $f(\bfx) = \bfdt_{[0,255]^n} (\bfx)$, $g(\bfu) = \| \bfu - \bfc \|^2/2$, $A = \R$, and $B = \D$.
	
	We conduct reconstruction experiment on four CT images: Shepp Logan \cite{shepp1974fourier}, FORBILD \cite{yu2012simulation}, and two samples selected from Brain CT Images with Intracranial Hemorrhage Masks dataset\footnote{https://www.kaggle.com/datasets/vbookshelf/computed-tomography-ct-images} (denoted as ICH1 and ICH2). The projection angles were chosen as $s$ equidistant samples over $(0,180^\circ]$, and then the sampled Radon transform by the build-in Matlab function $\texttt{radon}$. The information of data is summarized in Tab. \ref{tab-ct-dataset}.
	\begin{table}[htbp]
		\caption{Information of data for CT restoration}
		\label{tab-ct-dataset}
		\centering
		\begin{tabular}{cccc}
			\toprule
			Dataset & Image size $n$ & Total number of measurements $m$ & Number of angles $s$ \\
			\midrule
			Shepp Logan & $128 \times 128$ & 3700  & 20 \\
			FORBILD      & $200 \times 200$ & 5740  & 20 \\
			ICH1          & $200 \times 200$ & 14350 & 50 \\
			ICH2          & $200 \times 200$ & 14350 & 50 \\
			\bottomrule
		\end{tabular}
	\end{table}
	
	For the parameter setting of PPSNA, we take $\sigma = 10$ at the warm-start stage. After obtain a initial point $\bfom^l$, we set $\bfmu = \beta_1 \texttt{abs}(\bfv^l) + \beta_2$, where $\beta_1 = 100$ and $\beta_2 =  \min \{|(B\bfx)_i| /(mn) \ | \ (B\bfx)_i \neq 0 \}$. Other parameter setting is the same as that in Section \ref{exp-jump-sparse}. Pottslab, CVX, and PD3O are used for comparison in this experiment. For CT image restoration problem, Pottslab solves a model similar to \eqref{ct-restoration} without the box constraint \cite{storath2015joint}. CVX and PD3O solve convex optimization, which replaces $\pla$ in \eqref{ct-restoration} with $\ell_1$ loss with weight vector $\bfla$. Particularly, when reconstructing FORBILD, ICH1, and ICH2, the default solver SDPT3 of CVX runs out of memory of our laptop. Alternatively, we adopt solver SeDuMi for these three images. In the following experiments, we record PSNR and Time of the algorithms to evaluate their performance. We stop the algorithm if its Time is more than 7.5 hours. In this experiment, we take $\bfla = \lambda \bfone$. To achieve better performance of the algorithms, the regularization parameter $\lambda$ is selected from $\{ 10^{-3}, 2\times 10^{-3}, \cdots, 10^{-2}, 2 \times 10^{-2}, \cdots, 1, 2, \cdots, 10 \}$.
	
\begin{figure}[htbp]
		\centering
		\includegraphics[width=0.8\textwidth]{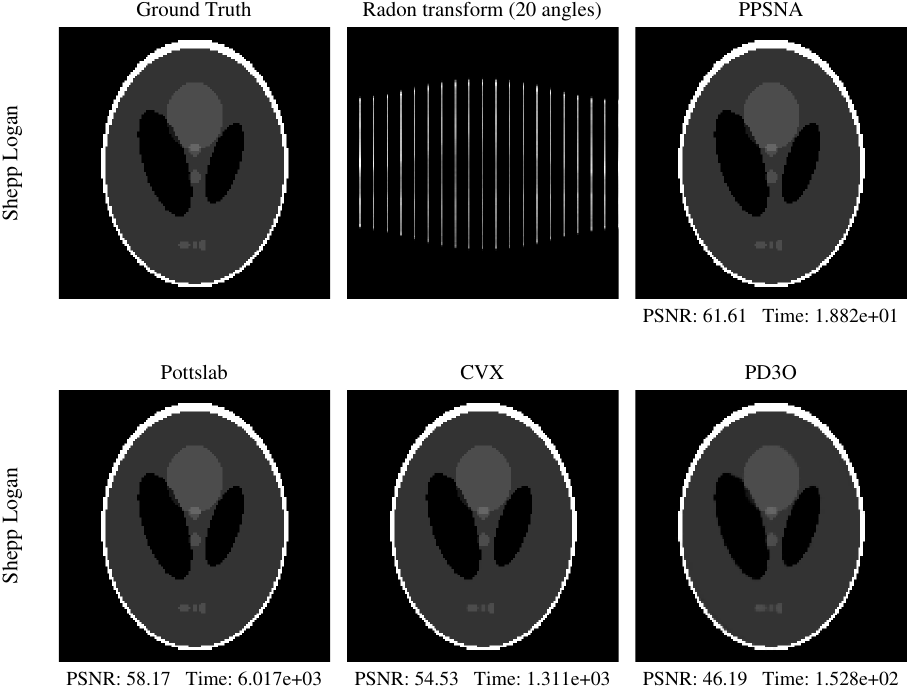}
		\caption{PSNR and Time of all algorithms on Shepp Logan.}
		\label{fig-shepp-logan}
	\end{figure}

\begin{figure}[htbp]
		\centering
		\includegraphics[width=0.8\textwidth]{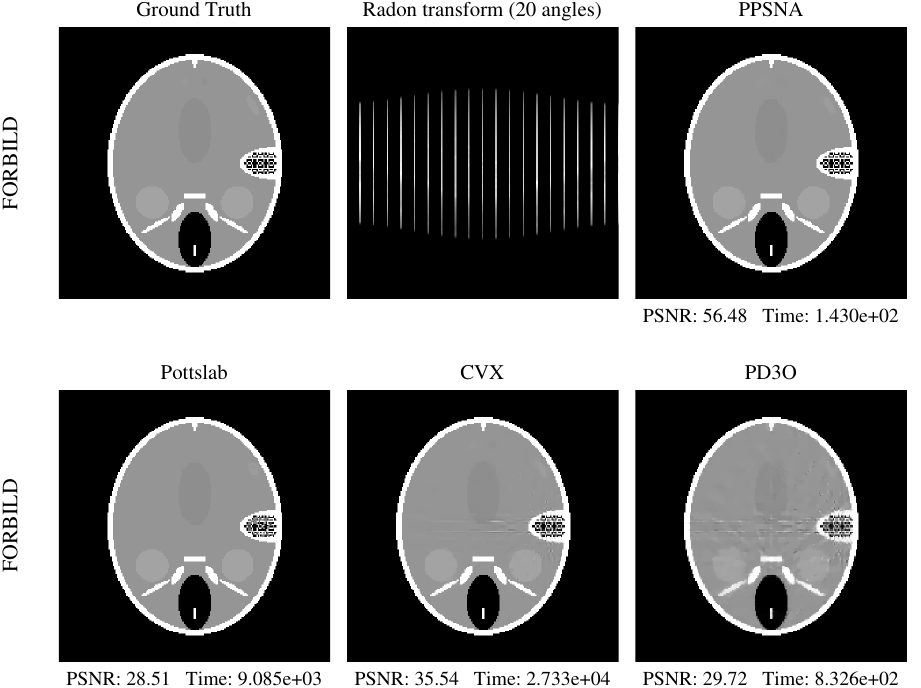}
		\caption{PSNR and Time of all algorithms on FORBILD.}
		\label{fig-forbild}
	\end{figure}

\begin{figure}[htbp]
		\centering
		\includegraphics[width=0.8\textwidth]{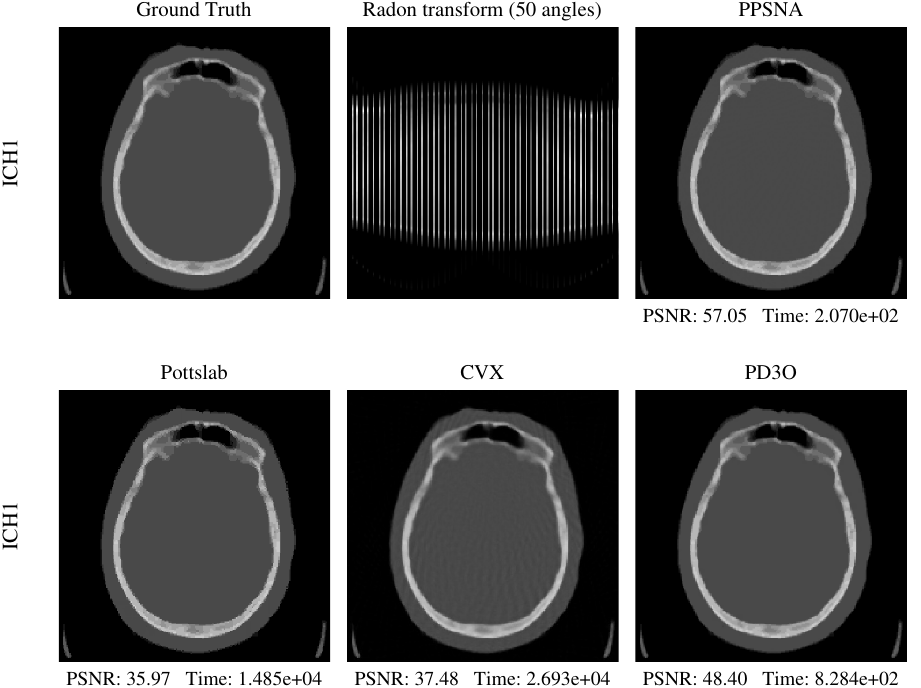}
		\caption{PSNR and Time of all algorithms on ICH1.}
		\label{fig-ich1}
	\end{figure}

\begin{figure}[htbp]
		\centering
		\includegraphics[width=0.8\textwidth]{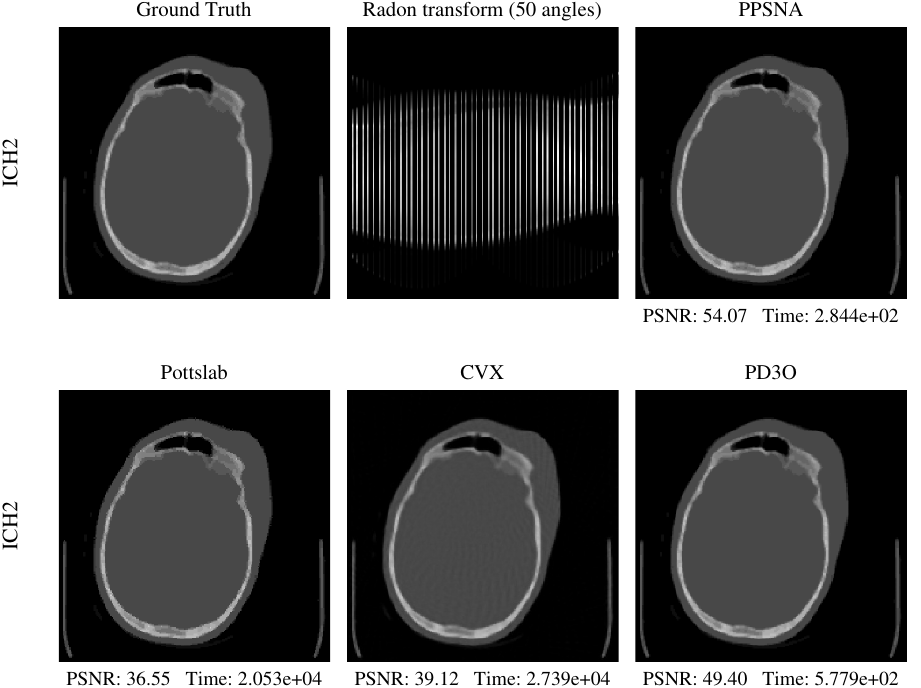}
		\caption{PSNR and Time of all algorithms on ICH2.}
		\label{fig-ich2}
	\end{figure}

\begin{figure}[htbp]
		\subfigure{
			\begin{minipage}[t]{0.47\linewidth}
				\centering
				\includegraphics[width=2.6in]{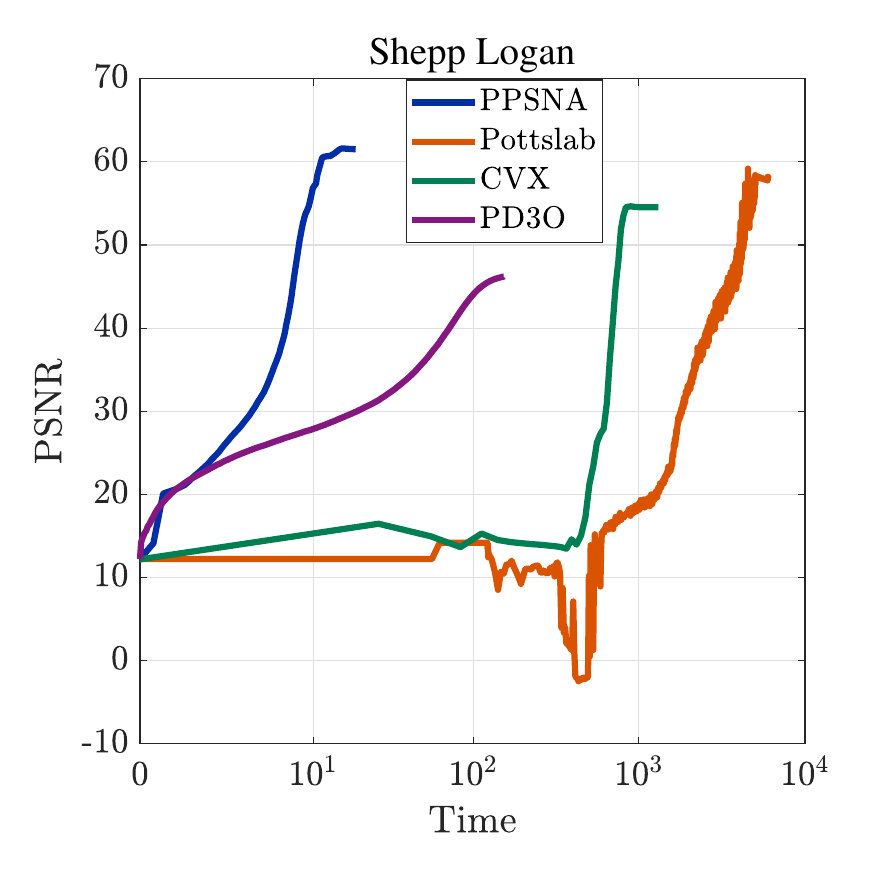}
			\end{minipage}%
			\begin{minipage}[t]{0.47\linewidth}
				\centering
				\includegraphics[width=2.4in]{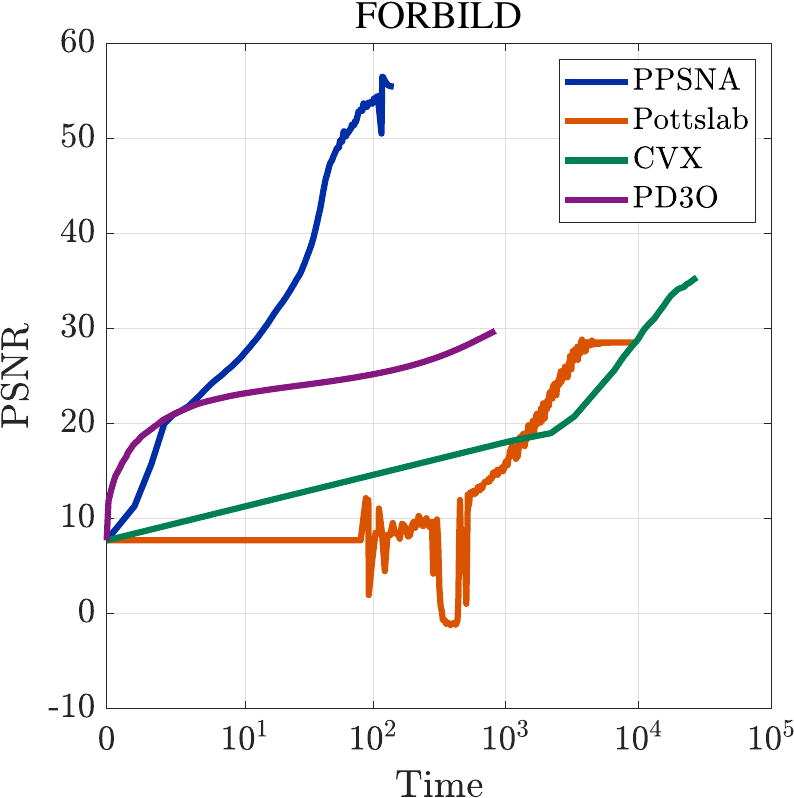}
			\end{minipage}%
		}
		
		\subfigure{
			\begin{minipage}[t]{0.47\linewidth}
				\centering
				\includegraphics[width=2.4in]{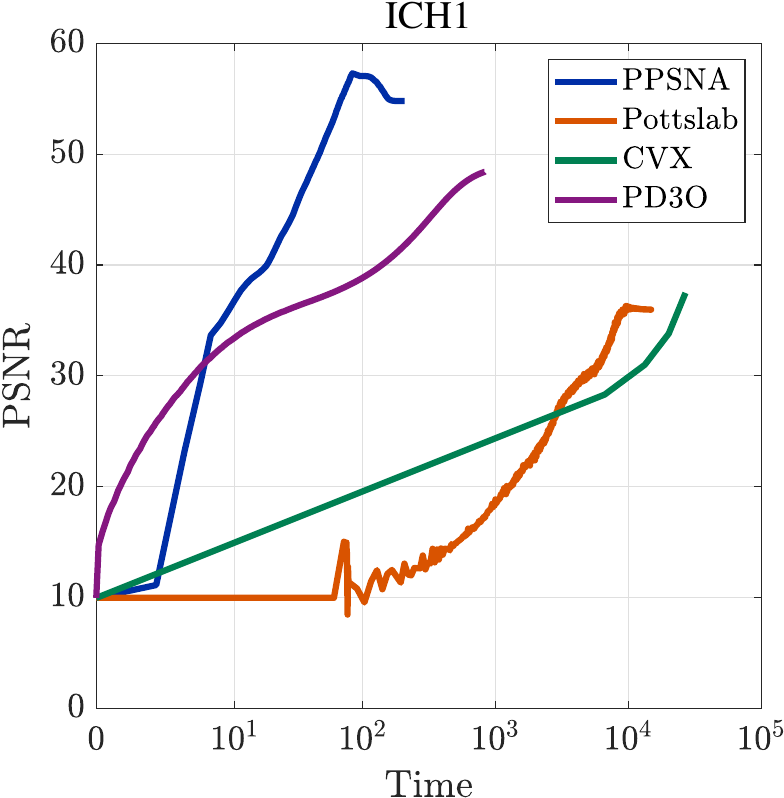}
			\end{minipage}%
			\begin{minipage}[t]{0.47\linewidth}
				\centering
				\includegraphics[width=2.4in]{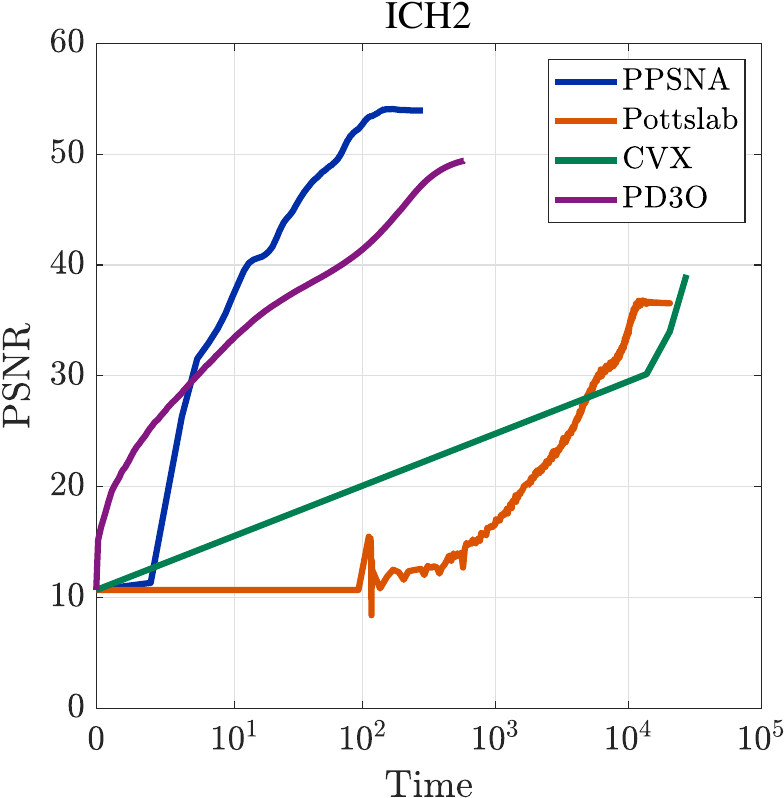}
			\end{minipage}%
		}
		\caption{The changes of PSNR along with Time for all the algorithms. 
			\label{psnr-time}}
		{}
	\end{figure}

	Finally, let us give some comments based on observations from Figs. \ref{fig-shepp-logan}-\ref{psnr-time}. For the dataset Shepp Logan, PPSNA and Pottslab achieve higher PSNR than the other two algorithms. Moreover, Time of PPSNA is much shorter than Pottslab and CVX. For the dataset FORBILD, PPSNA has the best reconstruction performance. Its Time less than 20\% of PD3O. For datasets ICH1 and ICH2, PPSNA and PD3O have better PSNR compared with the other two algorithms. Time of PPSNA is smaller than half of Pottslab. Overall, PPSNA and PD3O demonstrate fast increase on PSNR along with Time. In comparison, the PSNR of PD3O increases faster than that of PPSNA during the first 10 seconds, but more slowly thereafter.

	%
	
	\section{Conclusions} \label{sec-conclu}
	
	Composite Cardinality Optimization (CCOP) has been employed in a wide range of image processing problems. 
	Contrast to existing approaches such as relaxation, smoothing, and/or splitting methods, this paper addresses CCOP directly and proposes a general numerical framework based on its stationary dual problem. 		
	To design an efficient inner subproblem solver, we first use a coordinate descent step to identify a low-dimensional subspace, and then solve a convex program \eqref{Pk} by the classic semismooth Newton algorithm (SNA). To terminate the SNA early, we propose implementable inexact stopping criteria. We prove the global convergence of our inexact PPA under relative-interior-type condition, and then analyze the convergence rate of dual sequence under a local error bound condition. These conditions holds automatically in several common applications of image processing. This means we bypass the surjectivity-type assumptions widely used in global convergence of Lagrangian-based methods in nonconvex setting. Finally, the numerical experiment results on jump-sparse signal recovery and CT image restoration demonstrate the competitive performance of our proposed method.
	An interesting future research direction is to extend this framework to other models and methods \cite{hovhannisyan2016magma, tsipinakis2023multilevel}. In particular, if $\pla$ is replaced by a group-cardinality function, problem \eqref{P} can cover image-processing applications involving the isotropic $\ell_0$-TV regularizer. This extension will be investigated in our future work.

\section*{Acknowledgment}
This paper was supported by Hong Kong RGC General Research Fund (PolyU/15309223, 15303124), PolyU AMA Projects (P0044200, P0045347), National Key R\&D Program of China (2023YFA1011100), and 111 Project of China (B16002).

\section*{Appendix A. Proof of \cref{thm-SNA-global-convergence}}
The line search step of SNA ensures that $\{ \Psi_k(\obx^j) \}_{j \in \mbN}$ must be a decreasing sequence. Since $\Psi_k$ is coercive by \cite[Lemma 2.24]{planiden2016strongly}, then $\{ \obx^j \}_{j \in \mbN}$ is bounded. The boundedness of $\{ \obw^j \}_{j \in \mbN}$ follows from smoothness of $g_k$, \eqref{gradient-fghk} and \eqref{dual-recovery}.

(i) From \cite[Theorem 3.4]{zhao2010newton}, each accumulation point of $\{\obx^j\}_{j \in \mbN}$ is a global minimizer of \eqref{Pk}. Then \eqref{Pk-convergence} can be further derived by the convexity and smoothness of $\Psi_k$.


(ii) Let us denote $\bfp^j = \nabla f_k(\obx^j)$. Then by \eqref{gradient-fghk}, we have
\begin{align*}
	\nabla \Psi_k( \obx^j  ) = \bfp^j + A^\top \oby^j + B^\top \obz^j = \bfp^j + Q^\top\obw^j.
\end{align*}
Using \eqref{fenchel-ineq}, \eqref{conjugate-subgradient}, and \eqref{dual-recovery}, we have
\begin{align*}
	&g_k(A\obx^j) + g^*_k(\oby^j) = \langle \oby^j, A \obx^j \rangle,~ h_k(B\obx^j) + h^*_k(\obz^j) = \langle \obz^j, B \obx^j \rangle \\
	& f_k(\obx^j) + f^*_k(\bfp^j) = \langle \obx^j, \bfp^j \rangle = \langle \obx^j, \nabla \Psi_k (\bfx^{j}) - Q^\top \obw^j \rangle \\
	& f_k(\obx^j) + f^*_k(-Q^\top \obw^{{j}}) \geq - \langle \obx^j, Q^\top \obw^j \rangle
\end{align*}
Then we can give estimation of duality gap by using the above formulas:
\begin{align*}
	0\leq G_k(\obx^j, \obw^j) = \Psi_k(\obx^j) + \Xi_k(\obw^j) = f^*_k(-Q^\top\obw^j) - f^*_k (\bfp^j) + \langle \obx^j, \nabla \Psi_k(\obx^j) \rangle
\end{align*}
By boundedness of $\{ \obx^j \}_{j \in \mbN}$ and \eqref{Pk-convergence}, we have $\lim_{j \to \infty} \langle \obx^j, \nabla \Psi(\obx^j) \rangle = 0$. Since $f_k$ is coercive, then $\dom f^* = \mbR^n$ by \cite[Theorem 11.8]{RockWets98}. Then the convex function $f^*$ is uniformly continuous on any compact subset of $\mbR^n$. Therefore, $\lim_{j \to \infty} f^*_k(-Q^\top\obw^j) - f^*_k (\bfp^j) = 0$ because $\lim_{j \to \infty} \| \bfp^j + Q^\top\obw^j \| = \lim_{j \to \infty} \|\nabla \Psi_k( \obx^j  )\| = 0$ and $\{ (\bfp^j, Q^\top\obw^j) \}_{j \in \mbN}$ is bounded. Overall, we can derive $\lim_{j \to \infty} G_k(\obx^j, \obw^j) = 0$. 

(iii) From (i) and (ii), we can obtain $\lim_{j \to \infty} \Xi_k(\obw^j) = - \min\Psi_k$. The second line of \eqref{optimality-k} implies $\min\Xi_k + \min\Psi_k = 0$ and therefore \eqref{lim-Phik} holds. Since $\Xi_k$ is $(1/\gamma_k)$-strongly convex, we denote $\whw^{k+1}$ as its unique global minimizer and it holds that $\Xi_k (\obw^j) - \min\Xi_k \geq \frac{1}{2\gamma_k} \| \obw^j - \whw^{k+1} \|^2$. Finally, by using \eqref{lim-Phik}, we have $\lim_{j \to \infty} \obw^j = \whw^{k+1}$.
\ep


\section*{Appendix B. Proof of \cref{thm-SNA-loc-superlinear} }Before giving the proof of \cref{thm-SNA-loc-superlinear}, we need the following lemma for preparation.
\begin{lemma} \label{lem-uniform-nonsing}
	Suppose that Assumption \ref{asm-ssm-reg} (ii) holds, then there exists $\delta, \vartheta > 0$ such that 
	\begin{align*}
		\bfla_{\min}(H) \geq \vartheta,~~\mbox{whenever}~~ H \in \whp^2 \Psi_k(\bfx)~\mbox{and}~\bfx \in \N(\whx^{k+1},\delta).
	\end{align*}
\end{lemma}
\bp
By \cite[Proposition 7.1.4]{facchinei2003finite}, we know that $\whp^2 \Psi_k(\bfx)$ is nonempty, convex, and compact. We denote
\begin{align} \label{xi>0}
	\vartheta := \frac{1}{2} \inf \{ \bfla_{\min} (\whH) ~ | ~ \whH \in \whp^2 \Psi_k(\whx^{k+1}) \}
\end{align}
Since $\bfla_{\min}(\cdot)$ is continuous, the infimum above is attainable, and therefore $\vartheta > 0$ by Assumption \ref{asm-ssm-reg} (ii) and convexity of $f_k$ and $g_k$. Moreover, by the upper semicontinuity of $\partial^2 f_k$ and $\partial^2 g_k$, we know that there exists $\delta > 0$ such that for any $H \in \whp^2 \Psi_k(\bfx)$ and $\bfx \in \N(\whx^{k+1},\delta)$, there exists $\whH \in \whp^2 \Psi_k (\whx^{k+1})$ such that $\| H - \whH \| \leq \vartheta$. Then by Weyl's Perturbation theorem \cite[Corollary III.2.6]{bhatia2013matrix} and \eqref{xi>0}, we have
\begin{align*}
	\bfla_{\min}(H) \geq \bfla_{\min}(\whH) - | \bfla_{\min}(H) - \bfla_{\min}(\whH)| \geq \bfla_{\min}(\whH) - \| H - \whH \| \geq 2\vartheta - \vartheta = \vartheta
\end{align*}
The proof is finished.
\ep

Now we are ready to give the proof of \cref{thm-SNA-loc-superlinear}.

\bp From \eqref{partial-whp} and Assumption \ref{asm-ssm-reg} (ii), we know that $\bfd^\top H \bfd > 0$ for any $H \in \partial^2 \Psi_k(\whx^{k+1})$ and $\bfd \neq 0$. Then by \cite[Proposition 7.4.12]{facchinei2003finite} and the convexity of $\Psi_k$, problem \eqref{Pk} admits the unique global minimizer $\whx^{k+1}$ and thus $\lim_{j \to \infty} \obx^j = \whx^{k+1}$ by Thm. \ref{thm-SNA-global-convergence}.

Let us take $H_j = R_j + A^\top U_j A + \gamma_k B_{T_k:}^\top B_{T_k:}$, where $R_j \in \partial^2 f_k(\obx^j)$ and $U_j \in \partial^2 g_k(A\obx^j)$. Since Assumption \ref{asm-ssm-reg} (i) implies the semismoothness of $\nabla f_k$ and $\nabla g_k$, we have 
\begin{align*}
	& \| \nabla f_k(\obx^j) - \nabla f_k(\whx^{k+1}) - R_j(\obx^j - \whx^{k+1}) \| = o(\| \obx^j - \whx^{k+1} \|) \\
	& \| \nabla g_k(A\obx^j) - \nabla g_k(A\whx^{k+1}) - U_j(A\obx^j - A\whx^{k+1}) \| = o(\| \obx^j - \whx^{k+1} \|)
\end{align*}
Then by the representation \eqref{equation-k} of $\Psi_k$, we can further derive
\begin{align*} 
	\| \nabla \Psi_k(\obx^j) - \nabla \Psi_k(\whx^{k+1}) - H_j(\obx^j - \whx^{k+1}) \| = o(\| \obx^j - \whx^{k+1} \|) 
\end{align*}
We denote $\bfr^{j} = (H_j + \theta_j I) \bfd^j + \nabla \Psi_k (\obx^j)$ and then as the framework in Alg. \ref{SNA}, we have
\begin{align} \label{r-o}
	\| \bfr^j \| \leq \tau_2 \| \nabla \Psi_k (\obx^j) \|^{1+\tau} = \tau_2 \| \nabla \Psi_k (\obx^j) - \nabla \Psi_k (\whx^{k+1}) \|^{1+\tau} \leq O(\| \obx^j - \whx^{k+1} \|^{1+\tau})
\end{align}
where the inequality follows from Lipschitz continuity of $\nabla \Psi_k$. Similarly, we can prove
\begin{align*}
	\theta_j \| \obx^j - \whx^{k+1} \| = \tau_1 \| \nabla \Psi_k (\obx^j) - \nabla \Psi_k (\whx^{k+1}) \| \| \obx^j - \whx^{k+1} \| \leq O(\| \obx^j - \whx^{k+1} \|^2)
\end{align*}
Based on the above three formulas and Lem. \ref{lem-uniform-nonsing}, we have the following estimation:
\begin{align} 
	&	\| \obx^j + \bfd^j - \whx^{k+1} \| \notag \\ = & \| \obx^j - \whx^{k+1} + (H_j + \theta_j I)^{-1} (\bfr^j - \nabla \Psi_k(\obx^j))\| \notag \\
	\leq & (1/\vartheta) \| (H_j + \theta_j I)(\obx^j - \whx^{k+1}) + \bfr^j - \nabla \Psi_k(\obx^j) + \nabla\Psi_k(\whx^{k+1}) \| \notag \\ 
	\leq& (1/\vartheta) (\| \nabla \Psi_k(\obx^j) - \nabla \Psi_k(\whx^{k+1}) - H_j(\obx^j - \whx^{k+1}) \| + \| \bfr^j \| + \theta_j \| \obx^j - \whx^{k+1} \| ) \notag \\ 
	=& o(\| \obx^j - \whx^{k+1} \|). \label{superlinear}
\end{align}
This implies
\begin{align*}
	\| \bfd^{j} \| \geq \| \obx^j - \whx^{k+1} \| - \| \obx^j + \bfd^j - \whx^{k+1} \|  = O(\| \obx^j - \whx^{k+1} \|). 
\end{align*}
Combining the above inequality with \eqref{r-o}, we have $\| \bfr^j \| \leq O(\|\bfd^j\|^{1+\tau})$. Then by using Lem. \ref{lem-uniform-nonsing}, we can obtain the following inequality when $j$ is sufficiently large 
\begin{align*}
	\langle \nabla \Psi_k(\obx^j), \bfd^j \rangle = \langle \bfr^j, \bfd^j \rangle - \langle (H_j + \theta_j I)\bfd^j, \bfd^j \rangle \leq \| \bfr^j \| \| \bfd^j \| - \vartheta \| \bfd^j \|^2 \leq -(\vartheta/2) \| \bfd^j \|^2.
\end{align*}
It follows from \cite[Theorem 3.2]{facchinei1995minimization} that $\alpha_j = 1$ when $j$ is sufficiently large and therefore $\obx^{j+1} = \obx^j + \bfd^j$. This together with \eqref{superlinear} imply \eqref{SNA-superlinear}.

Finally, let us consider the dual sequence $\{ \obz^j \}_{j \in \mbN}$. From \eqref{dual-recovery} and $\gamma_k$-smoothness of $g_k$ and $h_k$, we have
\begin{align*}
	\| \oby^j - \why^{k+1} \| \leq \gamma_k \|A\| \| \obx^j - \whx^{k+1} \| ~~\mbox{and}~~ \| \obz^j - \whz^{k+1} \| \leq \gamma_k \|B\| \| \obx^j - \whx^{k+1} \|,
\end{align*}
Then we can derive the R-superlinear convergence rate of $\{ \obw^j \}_{j\in \mbN}$ by \eqref{SNA-superlinear}.
\ep

\section*{Appendix C. Proof of \cref{thm-finite-termination}}
Let us recall that $\whw^{k+1}$ is the unique global minimizer of \eqref{Dk} and $\bfw^k$ is the $k$-th iterate of PPSNA. 

We first consider the case that $\bfw^k \neq \whw^{k+1}$. According to Thm. \ref{thm-SNA-global-convergence}, we have
\begin{align*}
	&\lim_{j \to \infty} \Xi_k(\obw^j) - \Xi_k(\wtw^j) \leq \lim_{j \to \infty} \Xi_k(\obw^j) - \min \Xi_k = 0, \\
	& \lim_{j \to \infty} G(\obx^j, \obw^j) = 0, \\
	& \lim_{j \to \infty} \| \obw^j - \bfw^k \|_M = \| \whw^{k+1} - \bfw^k \|_M > 0.
\end{align*}
Then \eqref{C1}, \eqref{C2}, and \eqref{C3} must hold for $(\bfx^{k+1}, \bfw^{k+1}) = (\obx^j, \obw^j)$ when $j$ is sufficiently large.

Next suppose that $\bfw^k = \whw^{k+1}$. Since $f$ is coercive by Assumption \ref{asm-coercive-f}, it holds that $\dom f^* = \mbR^n$ according to \cite[Theorem 11.8]{RockWets98}. Then by the optimality condition of \eqref{Dk} at $\bfw^k$, we have 
\begin{align*}
	0 \in - \left( \begin{array}{c}
		A \\
		B
	\end{array} \right) \partial f^*( -Q^\top \bfw^k ) + \left( \begin{array}{c}
		\partial g^*(\bfy^k) \\
		\V^k
	\end{array} \right), 
\end{align*}
where $\V^k := \{ \bfv \in \mbR^p \ | \ \bfv_{T_k} = 0 \}$. Since $\bfw^k = \whw^{k+1}$ and $\whz^{k+1}_{\OT_k} = 0$, then $\V^k \subseteq \partial \pmu(\bfz^k)$ from \eqref{subdiff-l0}. This means $\bfw^k$ is exactly a stationary point of \eqref{D}, and therefore PPSNA should have been terminated at $\bfw^k$. This contradiction indicates it holds that $\bfw^k = \whw^{k+1}$ and this completes the proof. \ep

\section*{Appendix D. Proof of \cref{lem-c1}} 
(i) From \eqref{dual-recovery}, $\bfz^{k+1}_{\OT_k} = 0$ and hence $\pmu(\bfz^{k+1}) \leq \pmu(\wbz^k)$. By the coordinate minimization \eqref{coordinate-min} and \eqref{C1}, we can obtain
\begin{equation} \label{suff-des1} 
	\begin{aligned}
		\Xi(\bfw^k) + \pmu(\bfz^k) = & \Xi_k(\bfw^k) + \pmu(\bfz^k) \geq \Xi_k(\wtw^k) + \pmu(\wbz^k)  \\ 
		\geq & \Xi_k(\bfw^{k+1}) + \pmu(\bfz^{k+1}) - \frac{1}{4\gamma_k} \| \bfw^{k+1} - \bfw^k \|^2_M \\
		\geq &  \Xi(\bfw^{k+1}) + \pmu(\bfz^{k+1}) + \frac{1}{4\gamma_k} \| \bfw^{k+1} - \bfw^k \|^2_M 
	\end{aligned}
\end{equation}
This implies the sufficient descent property of $\{ \Xi(\bfw^k) + \pmu(\bfz^k) \}_{k \in \mbN}$. The sequence is also bounded because the objective function of \eqref{D} is lsc and $\{ \bfw^k \}_{k \in \mbN}$ is bounded. Then $\{ \Xi(\bfw^k) + \pmu(\bfz^k) \}_{k \in \mbN}$ must be convergent and its sufficient descent property further implies $\lim_{k \to \infty} \| \bfw^{k+1} - \bfw^k \| =0$ due to $\gamma_k \leq \gamma_\infty$ and positive definiteness of $M$.

(ii) We have the following estimation
\begin{align} 
	\big| \pmu(\bfz^{k+1}) - \pmu(\bfz^{k})  \big| \leq & \big| (\Xi(\bfw^{k+1}) + \pmu(\bfz^{k+1})) - (\Xi(\bfw^{k}) + \pmu(\bfz^{k})) \big| \notag \\
	&+ \big| f^*(-Q^\top\bfw^{k+1}) - f^*(-Q^\top\bfw^{k}) \big| + \big| g^*(\bfy^{k+1}) - g^*(\bfy^k) \big|  \label{cardz-upper-bound}
\end{align}
Since $f^*$ is convex and $\dom f^* = \mbR^n$, Assumption \ref{asm-bound} (i) indicates that $f^*$ is Lipschitz continuous on a compact set including $\{ -Q^\top \bfw^{k} \}_{k \in \mbN}$. Assumption \ref{asm-bound} (ii) implies $g^*$ is uniformly continuous on a compact set including $\{ \bfy^{k} \}_{k \in \mbN}$. Then from  \eqref{cardz-upper-bound} and assertion (i), we can obtain $\lim_{k \to \infty } | \pmu(\bfz^{k+1}) - \pmu(\bfz^{k})  \big| = 0$. Since the value of $\pmu(\cdot)$ can only take finitely many values, then $\pmu(\bfz^k)$ must remain the same when $k$ is large enough. 
\ep

\section*{Appendix E. Proof of \cref{lem-c1-c2}} (i) As \eqref{Dk} is a strongly convex program, we have
\begin{align*}
	\Xi_k( \wtw^k ) - \Xi_k( \whw^{k+1} ) \geq \frac{1}{2\gamma_k} \| \wtw^k - \whw^{k+1} \|^2 ~\mbox{and}~ \Xi_k( \bfw^{k+1} ) - \Xi_k( \whw^{k+1} ) \geq \frac{1}{2\gamma_k} \| \bfw^{k+1} - \whw^{k+1} \|^2.
\end{align*}
Adding the above inequalities implies
\begin{align*}
	\Xi_k ( \wtw^k ) - \Xi_k(\bfw^{k+1}) + 2\big( \Xi_k( \bfw^{k+1} ) - \Xi_k(\whw^{k+1}) \big) \geq & \frac{1}{2\gamma_k} \| \wtw^k - \whw^{k+1} \|^2 + \frac{1}{2\gamma_k} \| \bfw^{k+1} - \whw^{k+1} \|^2 \\
	\geq & \frac{1}{4\gamma_k} \| \bfw^{k+1} - \wtw^k \|^2
\end{align*}
By \eqref{C2} and weak duality between \eqref{Pk} and \eqref{Dk}, we have
\begin{align*}
	\Xi_k( \bfw^{k+1} ) - \Xi_k(\whw^{k+1}) \leq \Xi_k(\bfw^{k+1}) + \Psi_k( \bfx^{k+1} ) \leq \veps^2_k/(2\gamma_k).
\end{align*}
Then the above two inequalities imply 
\begin{align} \label{wk+1-wtw-upper-bound}
	\| \bfw^{k+1} - \wtw^k \|^2 \leq 4\gamma_k (\Xi_k ( \wtw^k ) - \Xi_k(\bfw^{k+1})) + 4\veps^2_k.
\end{align}
We can notice that \eqref{suff-des1} and Lems. \ref{lem-c1} (i) and (ii) yield $\lim_{k \to \infty} \Xi_k ( \wtw^k ) - \Xi_k(\bfw^{k+1}) = 0$. This together with \eqref{wk+1-wtw-upper-bound} lead to $\lim_{k\to \infty} \| \bfw^{k+1} - \wtw^k \|^2 = 0$. Finally, by using $\lim_{k \to \infty} \| \bfw^{k+1} - \bfw^k \| =0$, we can obtain \eqref{wtw-wk}. 

(ii) Let us define univariate functions:
\begin{align}
	&\psi_{ki}(t) := f^*\big( -Q^\top [\bfy^k; \wz^k_1; \cdots;\wz^k_{i-1}; t; z^k_{i+1};\cdots;z^k_m] \big), \\
	&\varphi_{ki}(t) := \frac{1}{2\gamma_k} \big\| [\bfy^k; \wz^k_1; \cdots;\wz^k_{i-1}; t; z^k_{i+1};\cdots;z^k_m] - \bfw^k \big\|^2_M,
\end{align}
and therefore $\bfxi_{ki}(t) = \psi_{ki}(t) + \varphi_{ki}(t)$. Assumption \ref{asm-bound} (i) and assertion (i) imply the boundedness of $\{ (\bfw^k, \wtw^k) \}_{k \in \mbN}$. Since $f^*$ is convex and $\dom f^* = \mbR^n$, then $f^*$ is Lipschitz continuous on any compact set. Therefore, we can find constant $\sigma_1$ such that $\psi_{ki}$ is $\sigma_1$-Lipschitz continuous on interval $[-2|\wz^k_i|,2|\wz^k_i|]$ for any $i \in [p]$ and $k \in \mbN$. Moreover, we can observe that $\phi_{ki}$ is $\sigma_2$-smooth with $\sigma_2:= \| M \|/\gamma_0$. Now let us consider $i \in T_k$. By definition of $T_k$ and the representation \eqref{sol-coordinate-min}, $\wz^k_i = \argmin_{t} \bfxi_{ki}(t)$ and we have
\begin{align*}
	0 \in \partial \psi_{ki}( \wz^k_i ) + \nabla \varphi_{ki} (\wz^k_i)  ~\mbox{and}~  \bfxi_{ki}(0) - \bfxi_{ki}(\wz^k_i) \geq \underline{\mu},
\end{align*}
where $\underline{\mu}:= \min_{i \in [p]} \mu_i$. The first formula above further implies $| \nabla \varphi_{ki} (\wz^k_i) | \leq \sigma_1$ from \cite[Theorem 9.13]{RockWets98}. By using the local  $\sigma_1$-Lipschitz continuity of $\psi_{ki}$ and $\sigma_2$-smoothness of $\varphi_{ki}$, the second formula above imply
\begin{align*}
	\underline{\mu} \leq \psi_{ki}(0) - \psi_{ki}(\wz^k_i) + \varphi_{ki}(0) - \varphi_{ki}(\wz^k_i) \leq & \sigma_1 |\wz^k_i| - \nabla \varphi_{ki} (\wz^k_i)\wz^k_i + \frac{\sigma_2}{2} (\wz^k_i)^2 \\
	\leq & 2\sigma_1 |\wz^k_i| + \frac{\sigma_2}{2} (\wz^k_i)^2.
\end{align*}
By direct computation, we can verify $|\wz^k_i| \geq ( -2\sigma_1 + \sqrt{4\sigma_1^2 + 2\sigma_2\underline{\mu}} )/\sigma_2 : = \vartheta_1$.

(iii) First, let us show that $\I(\bfz^k) \supseteq T_k$ for sufficiently large $k$. If this is not true, then there exists infinite index set $\K \subseteq \mbN$ such that for any $k \in \K$, we can find $i_k \in T_k$ with $\bfz^k_{i_k} = 0$. Then assertion (ii) implies $| \wz^k_{i_k} - z^k_{i_k} | \geq \vartheta_1$. By taking subsequence $\K_1 \subseteq \K$, we can find a fixed index $\overline{i}$ such that $| \wz^k_{\overline{i}} - z^k_{\overline{i}} | \geq \vartheta_1$ for all $k \in \K_1$. Therefore, we arrive at a contraction due to $\lim_{k \to \infty} \| \wtw^k - \bfw^k  \| = 0$. Now we claim that $\I(\bfz^k) \supseteq T_k \supseteq \I(\bfz^{k+1})$ when $k$ is large enough by \eqref{dual-recovery}. Finally, we can derive $\I(\bfz^k) = T_k = \I(\bfz^{k+1})$ by using Thm. \ref{lem-c1} (ii). \ep

\section*{Appendix F. Proof of \cref{thm-global-convergence}} 
(i) As we have stated that Alg. \ref{PPSNA} becomes a preconditioned PPA for \eqref{D-T-infty} when $k$ is sufficiently large. Then it follows from \cite[Theorem 2.3]{li2020asymptotically} that $\{\bfw^{k}\}_{k \in \mbN}$ converges to a global minimizer $\bfw^*$ of \eqref{D-T-infty} and therefore $0 \in \S(\bfw^*)$. Since $\V_\infty \subseteq  \partial \pmu(\bfz^*)$ by \eqref{subdiff-l0}, we can conclude that $\bfw^*$ is a stationary point of \eqref{D} by \eqref{stationary-D}. Moreover, it follows from Lem. \ref{lem-correspondence} that $\bfw^*$ is a local minimizer of \eqref{D}. 
	
	(ii) The $(k+1)$-th dual iterate updated by \eqref{dual-recovery} is as follows:
	\begin{align} \label{dual-recovery-k+1}
		\bfy^{k+1}  = \prox_{\gamma_k g^*} (\gamma_k A\bfx^{k+1} + \bfy^k),~~ \bfz^{k+1}_{T_k}  = \bfz^k_{T_k} + \gamma_k B_{T_k:} \bfx^{k+1} ~~ \mbox{and}~~ \bfz^{k+1}_{\OT_k} = 0. 
	\end{align}
	By \eqref{conjugate-subgradient} and \eqref{gradient-fghk}, we can obtain
	\begin{align*}
		g_k(A\bfx^{k+1}) + g^*_k(\bfy^{k+1}) = \langle \bfy^{k+1}, A \bfx^{k+1} \rangle,~ h_k(B\bfx^{k+1}) + h^*_k(\bfz^{k+1}) = \langle \bfz^{k+1}, B \bfx^{k+1} \rangle
	\end{align*}
	Denoting $\wtx^{k+1} := \prox_{\tfrac{1}{\gamma_k} f}(\bfx^{k+1} - Q^\top \bfw^k/ \gamma_k)$, then definition of the proximal operator implies
	\begin{align*}
		f_k(\bfx^{k+1}) = f(\wtx^{k+1}) + \frac{\gamma_k}{2} \| \wtx^{k+1} - ( \bfx^{k+1} - Q^\top \bfw^k/\gamma_k  ) \|^2 - \frac{1}{2\gamma_k} \| Q^\top \bfw^k \|^2.
	\end{align*}
	Next by using \eqref{C2} and the above equations, we have
	\begin{align*}
		\veps_k^2/(2\gamma_k) \geq& \Psi_k(\bfx^{k+1}) + \Xi_k(\bfw^{k+1}) = f_k(\bfx^{k+1}) + f^*_k(-Q^\top \bfw^{k+1}) + \langle \bfx^{k+1}, Q^\top\bfw^{k+1} \rangle \\
		= & f(\wtx^{k+1}) + \frac{\gamma_k}{2} \| \wtx^{k+1} - ( \bfx^{k+1} - Q^\top \bfw^k/\gamma_k  ) \|^2 - \frac{1}{2\gamma_k} \| Q^\top \bfw^k \|^2  + \langle \bfx^{k+1}, Q^\top\bfw^{k+1} \rangle \\
		& + f^*(-Q^\top \bfw^{k+1}) + \frac{1}{2\gamma_k} \| Q^\top (\bfw^{k+1} - \bfw^k) \|^2 \\
		= & f(\wtx^{k+1}) + f^*(-Q^\top \bfw^{k+1}) + \langle Q^\top \bfw^{k+1}, \wtx^{k+1} \rangle  \\
		& + \frac{1}{2\gamma_k} \|  Q^\top(\bfw^{k+1} - \bfw^k) - \gamma_k(\wtx^{k+1} - \bfx^{k+1}) \|^2 \geq 0, 
	\end{align*}
	where the last inequality follows from \eqref{fenchel-ineq}. Since $\lim_{k \to \infty} \bfw^k = \bfw^*$ and $f(\wtx^{k+1}) + f^*(-Q^\top \bfw^{k+1}) + \langle Q^\top \bfw^{k+1}, \wtx^{k+1} \rangle \geq 0$, passing $k$ to $\infty$ for the above formulas implies
	\begin{align} \label{lim-wtx}
		\begin{aligned}
			& \lim_{k \to \infty} f(\wtx^{k+1}) + f^*(-Q^\top \bfw^{k+1}) + \langle Q^\top \bfw^{k+1}, \wtx^{k+1} \rangle = 0\\ 
			& \lim_{k \to \infty} \| \wtx^{k+1} - \bfx^{k+1} \| = 0.
		\end{aligned}
	\end{align} 
	Now let us consider an accumulation point $\bfx^*$ of $\{ \bfx^{k} \}_{k \in \mbN}$. There exists infinite set $\K \subseteq \mbN$ such that $\lim_{k \in \K} \bfx^{k+1} = \bfx^* $. The second formula in \eqref{lim-wtx} indicates that $\lim_{k \in \K} \wtx^{k+1} = \bfx^* $. By lsc property of $f$ and $f^*$, the first formula in \eqref{lim-wtx} imply
	\begin{align*}
		0 \leq & f(\bfx^*) + f^*(-Q^\top \bfw^*) + \langle Q^\top \bfw^*, \bfx^* \rangle \\ \leq & \lim_{k \to \infty} f(\wtx^{k+1}) + f^*(-Q^\top \bfw^{k+1}) + \langle Q^\top \bfw^{k+1}, \wtx^{k+1} \rangle = 0,
	\end{align*}
	which means $-Q^\top \bfw^* \in \partial f(\bfx^*)$. Next we can take $k \to \infty$ for $k \in \K$ on each formula of \eqref{dual-recovery-k+1}. Then it follows from \cite[Theorem 1.25]{RockWets98} that
	\begin{align} 
		&\bfy^*  = \prox_{\gamma_\infty g^*} (\gamma_\infty A\bfx^* + \bfy^*), \label{prox-g*}\\
		&  B_{T_\infty:} \bfx^* = 0 ~~ \mbox{and}~~ \bfz^*_{\OT_\infty} = 0. \label{Bx-Tinfty=0}
	\end{align}
	By using \cite[Example 10.2]{RockWets98}, \eqref{prox-g*} implies $A\bfx^* \in \partial g^*(\bfy^*)$ and hence $\bfy^* \in \partial g(A\bfx^*)$ by \eqref{conjugate-subgradient}. Moreover, \eqref{Bx-Tinfty=0} indicates $\bfz^* \in \partial \pla(B\bfx^*)$ by \eqref{subdiff-l0}. Combining these two results with $-Q^\top \bfw^* \in \partial f(\bfx^*)$, we finally arrive at \eqref{stationary-P}. Finally, it follows from Lem. \ref{lem-correspondence} that $\bfx^*$ is a local minimizer of \eqref{P}. \ep
	
	\section*{Appendix G. Proof of \cref{thm-convergence-rate}} 
	As we have stated that Alg. \ref{PPSNA} becomes a preconditioned PPA for \eqref{D-T-infty} when $k$ is sufficiently large, and the global convergence of Alg. \ref{PPSNA} has been shown in Thm. \ref{thm-global-convergence}. Through the use of \cite[Theorme 2.5]{li2020asymptotically}, we can arrive at the desired conclusion.
	\ep
	
	\section*{Appendix H. Implementation of the coordinate minimization}
	\begin{proposition}
		The global minimizer of \eqref{coordinate-min} can be computed by \eqref{sol-coordinate-min}. 
	\end{proposition}
	\bp
	To compute the solution of \eqref{coordinate-min}, we only need to compare the optimal value of $\min_{z_i \in \mbR} \bfxi_{ki}(z_i) + \mu_i$ and $\bfxi_{ki}(0)$. The smaller one leads to the optimal value of \eqref{coordinate-min}. Since the two cases corresponds to solutions $\wtt^k_i := \argmin_{z_i \in \mbR} \bfxi_{ki}(z_i)$ and 0 respectively, we can obtain \eqref{sol-coordinate-min}.  
	\ep
	
	From \eqref{sol-coordinate-min}, we can observe that the main computational cost is from the one-dimensional minimization $\min \bfxi_{ki}$. In many cases, this step has closed-form expression. As for the two cases in our numerical experiment, we need to compute the one-dimension minimization with the following form:
	\begin{align} \label{quadratic-hinge-min}
		\argmin_{t \in \mbR} \bfxi(t) := at^2 + bt + \sum_{i=1}^n \max\{ c_i t - d_i, 0 \}, 
	\end{align}  
	where $a > 0$ and $b, d_i \in \mbR$ for $i \in [n]$. For convenience, we assume $c_i \neq 0$ for $i\in[n]$ because if $c_{i_0} = 0 $ for some $i_0 \in [n]$, then $\max\{ c_i t - d_i, 0 \} = \max\{ - d_i, 0 \}$ becomes a constant. Next let us give the closed form expression of \eqref{quadratic-hinge-min}.
	
	To proceed, we denote univariate functions $\phi_i(t):= \max\{ c_i t - d_i, 0 \}$ for $i \in [n]$. Then the non-differentiable point of $\phi_i$ is $\tau_i:= d_i/c_i$.
	After that, we sort $\{ \tau_i \}_{i\in[n]}$ in ascending order and select the distinct value by $\{ \theta_j \}_{j \in [r]}$, which means $\theta_1 < \theta_2 < \cdots < \theta_r$. In particular, the non-differentiable points of $\bfxi$ are $\{ \theta_j \}_{j \in [r]}$. For convenience, we further set $\theta_0:=-\infty$ and $\theta_{r+1}:=+\infty$. Let us denote
	\[
	q_0:=b+\sum_{\{ i:\ c_i<0 \}} c_i,
	\]
	Then for $t \in (\theta_0, \theta_1)$, we can compute $\partial \bfxi(t) = 2a + q_0$. For $j \in [r]$, define
	\[
	\Delta_j:=\sum_{\{ i:\, \tau_i=\theta_j \}} |c_i|,
	\]
	and then recursively we compute 
	\[
	q_j:=q_{j-1}+\Delta_j,~\mbox{for}~ j \in [r].
	\]
	For $t \in (\theta_j, \theta_{j+1})$ and $j = 0,1,\cdots,r$, we have
	\begin{align} \label{subdiff1-xi}
		\partial \bfxi (t) = 2at + q_j
	\end{align}
	Moreover, since $\bfxi$ is convex and it is also differentiable in $(\theta_{j-1}, \theta_{j})$ and $(\theta_j, \theta_{j+1})$ for $j \in [r]$. Then we can compute
	\begin{align} \label{subdiff2-xi}
		\partial \bfxi(t) = [2at+q_{j-1}, 2at+q_j]~~\mbox{for}~~j\in[r].
	\end{align}
	For convenience, we also denote $s_j := - q_j/(2a)$, which is strictly decreasing for $j \in \{ 0,1,\cdots,r \}$.
	Now we are ready to give the closed-form optimal solution of \eqref{quadratic-hinge-min}.
	\begin{proposition} \label{prop-sol-quadratic-hinge}
		Problem \eqref{quadratic-hinge-min} admits the unique optimal solution $t^*$, which must satisfy one of the following two cases.
		
		(i) If $s_j \in ( \theta_j, \theta_{j+1} )$ for some $j \in \{0,1,\cdots,r\} $, then $t^* = s_j$.
		
		(ii) If $\theta_j \in [s_{j}, s_{j-1}]$ for some $j \in \{1,2,\cdots,r\} $, then $t^*=\theta_j$.
	\end{proposition}
	\bp If $s_j \in ( \theta_j, \theta_{j+1} )$ holds for some $j \in \{0,1,\cdots,r\}$, then $s_j = -q_j/(2a)$ and \eqref{subdiff1-xi} imply
	\begin{align*}
		0 = 2as_j + q_j = \partial \bfxi(s_j).
	\end{align*} 
	This means $t^* = s_j$ is the unique global minimizer of \eqref{quadratic-hinge-min} by strong convexity of $\bfxi$.
	
	If $\theta_j \in [s_{j}, s_{j-1}]$ for some $j \in \{1,2,\cdots,r\} $, then $s_j = -q_j/(2a)$ and \eqref{subdiff2-xi} yields
	\begin{align*}
		0 \in [2a\theta_j + q_{j-1}, 2a\theta_j + q_{j}] = \partial \bfxi(\theta_j),
	\end{align*}
	which means $t^* = \theta_j$.
	
	Finally, it is noteworthy that $t^* \in (\theta_j, \theta_{j+1})$ for some $j \in \{0,1,\cdots,r\} $ in case (i) while $t^* = \theta_j$ for some $j \in \{1,2,\cdots,r\}$ in case (ii). Therefore, exactly one of the two cases holds. \ep
	
	At the end of the appendix, we discuss the computational complexity of solving \eqref{quadratic-hinge-min}. Computing $\{\tau_i\}_{i\in[n]}$ requires $O(n)$ operations. Sorting this sequence to obtain the distinct ordered values $\{\theta_j\}_{j\in[r]}$ costs $O(n\log n)$ operations, and computing $\{s_j\}_{j\in[r]}$ requires at most $O(n)$ additional operations. Hence, the overall complexity is $O(n\log n)$. It is noteworthy that, in \eqref{quadratic-hinge-min}, we assume $c_i\neq 0$ for all $i\in[n]$. Therefore, if $B$ is sparse, only the nonzero components need to be considered, and the effective value of $n$ can be substantially reduced. For instance, when $B=D_n$ (see \eqref{Dn}), each row of $B$ contains only two nonzero entries. In this case, the effective dimension is much smaller, and thus the computational complexity of solving \eqref{quadratic-hinge-min} is $O(1)$.

\bibliographystyle{siamplain}
\bibliography{CCOPT_ref} 

\end{document}